\documentclass[a4j,10pt]{article}
\usepackage{amsmath,amscd,amssymb}

\newcommand{\nbiga}{\mathcal{A}}

\newcommand{\nbigc}{\mathcal{C}}

\newcommand{\nbige}{\mathcal{E}}
\newcommand{\nbigf}{\mathcal{F}}
\newcommand{\nbigg}{\mathcal{G}}
\newcommand{\nbigh}{\mathcal{H}}

\newcommand{\nbigk}{\mathcal{K}}
\newcommand{\nbigl}{\mathcal{L}}
\newcommand{\nbigm}{\mathcal{M}}

\newcommand{\nbigo}{\mathcal{O}}
\newcommand{\nbigp}{\mathcal{P}}

\newcommand{\nbigr}{\mathcal{R}}
\newcommand{\nbigs}{\mathcal{S}}
\newcommand{\nbigt}{\mathcal{T}}
\newcommand{\nbigu}{\mathcal{U}}

\newcommand{\seisuu}{\Bbb{Z}}

\newcommand{\cnum}{{\bf C}}

\newcommand{\hyperh}{\Bbb{H}}

\newcommand{\gminip}{\frak p}
\newcommand{\gminiq}{\frak q}

\newcommand{\vece}{{\boldsymbol e}}

\newcommand{\vecv}{{\boldsymbol v}}
\newcommand{\vecu}{{\boldsymbol u}}

\newcommand{\veca}{{\boldsymbol a}}

\newcommand{\rarr}{\rightarrow}
\newcommand{\lrarr}{\longrightarrow}

\newcommand{\pf}{{\bf Proof}\hspace{.1in}}
\newcommand{\qed}{\mbox{\rule{1.2mm}{3mm}}}

\def\rank{\mathop{\rm rank}\nolimits}
\def\Spec{\mathop{\rm Spec}\nolimits}

\def\ker{\mathop{\rm Ker}\nolimits}
\def\length{\mathop{\rm length}\nolimits}

\def\Sym{\mathop{\rm Sym}\nolimits}

\def\Supp{\mathop{\rm Supp}\nolimits}

\newcommand{\nhom}{{\mathcal Hom}}

\newcommand{\nN}{\nbigp}
\newcommand{\mM}{\nbigm}
\newcommand{\nnN}{\nbigp_0}
\newcommand{\mmM}{\nbigm_0}

\newcommand{\TT}{\nbigu}
\newcommand{\TTalpha}{\TT_{\alpha}}
\newcommand{\LLA}{\seisuu}
\newcommand{\LLB}{C}
\newcommand{\TTt}{\TT_{\geq\,0}}

\newcommand{\bdmath}{\begin{displaymath}}
\newcommand{\edmath}{\end{displaymath}}

\newcommand{\beqn}{\begin{equation}}
\newcommand{\eeqn}{\end{equation}}

\newcommand{\beqnarray}{\begin{eqnarray}}
\newcommand{\eeqnarray}{\end{eqnarray}}

\newcommand{\bitemize}{\begin{itemize}}
\newcommand{\eitemize}{\end{itemize}}

\newcommand{\benumerate}{\begin{enumerate}}
\newcommand{\eenumerate}{\end{enumerate}}

\newcommand{\bdescriprion}{\begin{description}}
\newcommand{\edescriprion}{\end{description}}

\newtheorem{thm}{Theorem}[section]
\newtheorem{cor}{Corollary}[section]

\newtheorem{rem}{Remark}[section]
\newtheorem{lem}{Lemma}[section]
\newtheorem{prop}{Proposition}[section]

\def\longuparrow{\vcenter{%
   \lineskip0pt\lineskiplimit0pt\baselineskip0pt\ialign{%
   \hfil{##}\hfil\crcr\hbox to 0pt{\hss$\uparrow$\hss}\cr%
   \hbox to 0pt{\hss\vrule width.4pt depth 0pt height 1em\hss}\cr}}}
\def\longdownarrow{\vcenter{%
   \lineskip0pt\lineskiplimit0pt\baselineskip0pt\ialign{\hfil{##}\hfil%
   \crcr\hbox to 0pt{\hss\vrule width.4pt depth 0pt height 1em\hss}\cr%
   \hbox to 0pt{\hss$\downarrow$\hss}\cr}}}

\begin{document}

\title{The geometry of the parabolic Hilbert schemes\\
}
\author{Takuro Mochizuki\footnote{
takuro@sci.osaka-cu.ac.jp,
takuro@math.ias.edu}
}
\date{}

\maketitle

\abstract{
Let $X$ be a smooth projective surface and $D$ be a smooth divisor
over an algebraically closed field $k$.
In this paper, we discuss the moduli schemes
of the ideals of points of $X$
with parabolic structures at $D$.
They are called parabolic Hilbert schemes.
The first result is that the parabolic Hilbert schemes are smooth.
And then some of the studies of Ellingsrud-Str\o{}mme,
G\"ottsche, 
Cheah, Nakajima and Grojnowski
on the Hilbert schemes
can be naturally generalized in the parabolic case.

We determine the class of the parabolic Hilbert schemes
in the Grothendieck ring of $k$-varieties:
The class is described in terms of 
products of the symmetric powers
of $X$ and $D$
and the affine spaces.
Thus we obtain a formula for the generating functions
of the $E$-polynomials or the Poincar\'e polynomials
of the parabolic Hilbert schemes
of a smooth projective surface $X$ and a divisor $D$
over the complex number field $\cnum$.

Moreover we obtain the extension of the Nakajima-Grojnowski theory
for the parabolic Hilbert schemes:
i.e.,
we introduce the incidence varieties which induce the  operations
on the cohomology groups of the parabolic Hilbert schemes.
We see that the Heisenberg relations holds.
Thus we obtain a representation of Heisenberg algebra.

\noindent
Keywords: 
Hilbert scheme,
parabolic sheaf,
generating function.

\noindent
MSC: 14F45, 14M99.
}


\section{Introduction}

\subsection{Preliminary}
\subsubsection{Parabolic Hilbert schemes and our purpose in this paper}

Let $X$ be a smooth projective surface over an algebraically closed
field $k$,
and $D$ be a smooth divisor of $X$.
Let $E$ be a torsion-free sheaf.
We consider a filtration of $\nbigo_D$-modules
$\nbigg_{\ast}=\{\nbigg_{\alpha}\,|\,\alpha\in\seisuu\}$,
where $\nbigg_{\alpha}$ is a quotient sheaf of $\nbigg_{\alpha-1}$
and there are numbers $\alpha_-$ and $\alpha_+$ satisfying the following:
$\nbigg_{\alpha}=E\otimes\nbigo_D$
for any $\alpha<\alpha_-$
and $\nbigg_{\alpha}=0$ for any $\alpha\geq \alpha_+$.
Such filtration is called a parabolic structure.
A tuple $(I,\nbigg_{\ast})$
of ideal $I$ of points on $X$ and a parabolic structure
$\nbigg_{\ast}$ of $I$ at $D$
is called a parabolic ideal of points on $(X,D)$.
The moduli schemes of parabolic ideals of points
are called the parabolic Hilbert scheme of points on $X$.

We put $\nbigk_{\alpha}:=\ker(\nbigg_{\alpha}\lrarr\nbigg_{\alpha+1})$
and $\nbigo_Z:=\nbigo/I$
for a parabolic ideal $(I,\nbigg_{\ast})$.
We have the unique jumping number $\alpha_0$
such that $\nbigk_{\alpha_0}$ is not torsion
as an $\nbigo_D$-module.
In this paper,
we consider only the parabolic ideals whose jumping numbers are $0$.
Let $O$ be a point of $X$.
The condition
$\Supp(\nbigo_X/I)\cup
 \bigcup_{\alpha\neq 0}\Supp(\nbigk_{\alpha})\subset\{O\}$
determines the closed subschemes of the parabolic Hilbert schemes.
Such closed subschemes are called the punctual parabolic Hilbert schemes.

The starting point of the study is the following:
\begin{thm}[Corollary \ref{cor;smooth}]\label{thm;12.2.5}
Let $X$ be a smooth surface and $D$ be a smooth divisor.
The parabolic Hilbert schemes of points on $(X,D)$
are smooth.
\end{thm}

Thus we obtain a family of smooth varieties
which are moduli scheme of some algebro-geometric objects.
It seems natural to expect some interesting properties
for such a family as in the case of usual Hilbert schemes
of points on a smooth surface.
The study of Hilbert schemes of points on a smooth projective surface
is one of the most exciting subjects in the recent mathematics.
What the author would like to do first
in the study of parabolic Hilbert schemes
is to generalize the excellent researches for the Hilbert schemes
of points on a smooth surface.
In this paper,
we will discuss the following:
\begin{itemize}
\item
 We see the structure of the punctual parabolic Hilbert schemes.
 We obtain the cell decomposition and an extension of
 the theorem of Brian\c{c}on.
\item
 We consider the classes of the parabolic Hilbert schemes
 in the Grothendieck ring $K_0(V_k)$ of 
 the smooth varieties over $k$.
 We will obtain the formula to describe
 the parabolic Hilbert schemes in terms of the symmetric powers
 of $X$ and $D$.
 In particular, we obtain the formula of the Betti number
 of the parabolic Hilbert schemes.
 It is a generalization of the formula of G\"ottsche. 
 \item
 We see the operators obtained as the correspondence map
 induced by some incidence varieties.
 We will obtain the extension of Nakajima-Grojnowski theory
 in our parabolic Hilbert schemes.
\end{itemize} 

\subsubsection{Some notation}

We denote the set of integers by $\seisuu$.
We denote the set of integers larger (resp. less)
than $i$ by $\seisuu_{\geq\,i}$ (resp. $\seisuu_{\leq\,i}$).
We put $\TTalpha=\seisuu$ for any $\alpha\in \seisuu$.
We put as follows:
\[
 \TT:=\bigoplus_{\alpha\in\seisuu}\TTalpha,
\quad
 \TT_-:=\bigoplus_{\alpha<0}\TTalpha,
\quad
 \TT_+:=\bigoplus_{\alpha>0}\TTalpha.
\]
We have the natural direct sum decomposition
$\TT=\TT_0\oplus\TT_+\oplus \TT_-$.
We denote the projection of $\TT$ onto $\TTalpha$
by $\rho_{\alpha}$ for any integer $\alpha\in \seisuu$.
The restrictions of $\rho_{\alpha}$ $(\alpha\in\seisuu)$
to $\TT_+$ and $\TT_-$ are also denoted by $\rho_{\alpha}$.
We denote the projections of $\TT$ onto $\TT_+$
and $\TT_-$ by $\rho_+$ and $\rho_-$ respectively.
We put $\rho:=\rho_-\oplus \rho_+$.

For any element $\vecv\in \TT$,
we put as follows:
\[
 |\vecv|:=\sum_{\alpha\in\seisuu} |\rho_{\alpha}(\vecv)|,
\quad
 |\vecv|_+:=\sum_{\alpha>0}|\rho_{\alpha}(\vecv)|,
\quad
 |\vecv|_-:=\sum_{\alpha<0}|\rho_{\alpha}(\vecv)|,
\quad
d(\vecv):=2\rho_0(\vecv)+|\vecv|_+-|\vecv|_-.
\]
We put
$\nbiga:=
 \{\vecv\in\nbigt\,|\,\rho_{\alpha}(\vecv)\geq 0\,\,(\alpha\in\seisuu),
 \,\,\,
 \rho_0(\vecv)-|\vecv|_-\geq 0\}$.

Let 
$\vece_{\alpha}$ be the element of $\TT$
such that 
$\rho_{\beta}(\vece_{\alpha})$ is $0\,\,(\beta\neq\alpha)$
or $1\,\,(\beta=\alpha)$.
\label{subsubsection;7.6.10}
We put as follows:
\[
 C_{\alpha}:=
 \left\{
 \begin{array}{ll}
 \{m\cdot\vece_0+\vece_{\alpha}\,|\,m\in\seisuu_{\geq\,0}
  \}, &(\alpha>0), \\
  \{m\cdot\vece_0+\vece_{\alpha}\,|\,m\in\seisuu_{\geq\,1}
  \}, &(\alpha<0), \\
 \{m\cdot\vece_0\,|\,m\in\seisuu_{\geq\,1}
  \}, &(\alpha=0). \\
 \end{array}
 \right.
\]
We put $C_+:=\coprod_{\alpha>0}C_{\alpha}$,
$C_-:=\coprod_{\alpha<0}C_{\alpha}$
and $C=C_-\sqcup C_+\sqcup C_0$.
The set $-C$ is defined to be
$\{\vecv\in \TT\,|\,-\vecv\in C\}$.
Similarly we obtain the sets
$-C_0$, $-C_+$, and $-C_-$.
We put $\tilde{C}:=C\cup(-C)$.
Similarly the subsets
$\tilde{C}_+$, $\tilde{C}_0$ and $\tilde{C}_-$
are defined.

The function $g$ 
on the set $\tilde{C}$
is defined as follows:
\begin{equation}
 g(\vecu)=
 \left\{
   \begin{array}{rl}
 -1 & (\vecu\in -C_+,\mbox{ or }\vecu\in C_-)\\
 0  & (\mbox{otherwise}).\\
   \end{array}
 \right.
\end{equation}
For an element $\vecu$ of $\tilde{C}$,
let $\max(0,\vecu)$ denote $\vecu$ $(\vecu\in C)$
or $0$ $(\vecu\in -C)$.

Let $x_{\alpha}$ $(\alpha\in\seisuu)$ be variables.
For any finite subset $\Lambda\in \seisuu$,
we have the ring of formal power series
$\nbigr_{\Lambda}:=\seisuu[[x_{\lambda}\,(\lambda\in \Lambda)]]$.
For two finite subsets $\Lambda_i$ $(i=1,2)$ of $\seisuu$
such that $\Lambda_1\subset \Lambda_2$,
we have the natural morphism
$\nbigr_{\Lambda_2}\lrarr\nbigr_{\Lambda_1}$.
Thus we put $\nbigr:=\varprojlim_{\Lambda\subset\seisuu}\nbigr_{\Lambda}$.
Let $\vecv$ be an element of $\TT$ satisfying $\rho_{\alpha}(\vecv)\geq 0$.
The element $\prod_{\alpha}x_{\alpha}^{\rho_{\alpha}(\vecv)}$ of $\nbigr$
is denoted by $x^{\vecv}$.

\subsubsection{The length of parabolic ideals}

Let $x$ be a point of $X$ and $\nbigf$ be a coherent sheaf on $X$
with finite length.
We denote the length of $\nbigf$ at the point $x$
by $\length_x(\nbigf)$.
We denote the length of $\nbigf$ by $\length(\nbigf)$,
that is $\sum_{x\in X}\length_x(\nbigf)=\length(\nbigf)$.

Let $(I,\nbigg_{\ast})$ be a parabolic ideal of points
on $(X,D)$.
Let $x$ be a point of $X$.
Then we have the numbers
$\length_x(\nbigk_{\alpha})\in \TTalpha$ for any $\alpha\neq 0$
and $\length_x(\nbigo_Z)\in\TT_0$.
Thus we obtain the element of $\TT$.
We denote it by $\length_x(I,\nbigg_{\ast})$.
\begin{lem} \label{lem;7.4.1}
$\length_x(I,\nbigg_{\ast})$ is contained 
in $\nbiga$.
\end{lem}
\pf
We put $K:=\ker(I\otimes\nbigo_D\lrarr\nbigg_0)$.
We have to show the inequality
$\length_x(K)\leq \length_x(\nbigo_Z)$.
We denote the torsion part of $I\otimes\nbigo_D$
by $T$.
Then we only have to show the inequality
$\length_x(T)\leq \length_x(\nbigo_Z)$.
The sheaf $T$ is isomorphic to $Tor^1(\nbigo_Z,\nbigo_D)$.
We have an injection
$Tor^1(\nbigo_Z,\nbigo_D)\lrarr \nbigo_Z\otimes\nbigo(-D)$.
Thus we obtain the inequality desired.
\hfill\qed

We put 
$\length(I,\nbigg_{\ast}):=\sum_{x\in X}\length_x(I,\nbigg_{\ast})$.
We call it the length of parabolic ideal $(I,\nbigg_{\ast})$.
It is also the element of $\nbiga$.

\subsubsection{The parabolic Hilbert schemes and the symmetric powers
 associated with $\vecv$}

In general, we denote the $n$-th symmetric power of a scheme $Y$
by $Y^{(n)}$ for any integer $n$.
We have the natural projection $Y^n\lrarr Y^{(n)}$.
We denote the image of $(x_1,\ldots,x_n)\in Y^n$ via the projection
by $\sum_i x_i$.
It can be rewritten as
$\sum_{x\in X}f(x)\cdot x$.

The variety $Y_{\alpha}$ is defined to be
$X$ $(\alpha=0)$       
or
$D$ $(\alpha\neq 0)$.
Let $\vecv$ be an element of $\nbigt$
satisfying $\rho_{\alpha}(\vecv)\geq 0$ for any $\alpha$.
Then we put 
$X^{(\vecv)}
 :=\prod_{\alpha\in\seisuu} Y_{\alpha}^{(\rho_{\alpha}(\vecv))}$.
On the other hand,
we denote the moduli scheme of parabolic ideals
of points with parabolic length $\vecv$
by $X^{[\vecv]}$
for any element $\vecv\in \nbiga$.
In particular, if $\vecv=n\cdot\vece_0$,
then we denote $X^{[\vecv]}$ also by $X^{[n]}$.

In general, we put
$\Supp(\nbigf):=\sum_{x\in X}\length_x(\nbigf)\cdot x
 \in X^{(\length(\nbigf))}$
for a coherent sheaf $\nbigf$ on $X$.
For a parabolic ideal $(I,\nbigg_{\ast})$,
we have the points
$\Supp(\nbigo_Z)\in X^{(\length(\nbigo_Z))}$,
and $\Supp(\nbigk_{\alpha})\in D^{(\length(\nbigk_{\alpha}))}$
$(\alpha\neq 0)$.
Thus we have the natural morphisms
$X^{[\vecv]}\lrarr Y_0^{(\rho_0(\vecv))}$
given by the correspondence
$(I,\nbigg_{\ast})\longmapsto \Supp(\nbigo_Z)$,
and $X^{[\vecv]}\lrarr Y_{\alpha}^{(\rho_{\alpha}(\vecv))}$
given by the correspondence
$(I,\nbigg_{\ast})\longmapsto \Supp(\nbigk_{\alpha})$.
Thus we have the Hilbert-Chow morphism
$X^{[\vecv]}\lrarr X^{(\vecv)}$.


\subsection{The main results}

\subsubsection{Poincar\'e Polynomials}

One of our purpose in this paper is to obtain the formula
to describe the class of the parabolic Hilbert schemes
in terms of the symmetric powers of $X$ and $D$,
in the Grotendieck group $K_0(V_k)$ of smooth $k$-varieties
over $k$.
Let $A^{l}$ denote an $l$-dimensional affine space over $k$.
For a variety $Y$ and an integer $a$,
we put as follows:
\[
 S_a(Y,t)=\sum_{n=0}^{\infty} (Y\times A^a)^{(n)}\cdot t^n
 \in K_0(V_k)[[t]].
\]
\begin{thm}[G\"ottsche \cite{g2}]
Let $X$ be a smooth surface over an algebraically closed field $k$.
In $K_0(V_k)[[t]]$,
the following equality holds:
\begin{equation} \label{eq;7.2.1}
 \sum_{n=0}^{\infty}X^{[n]}\cdot t^n=
 \prod_{i=1}^{\infty}
 S_{i-1}(X,t^{i}).
\end{equation}
When the ground field $k$ is the complex number field $\cnum$,
we can obtain the formula of the generating function
of the $E$-polynomials of $X^{[n]}$
in terms of those of the symmetric powers of $X$.
\hfill\qed
\end{thm}
\begin{rem}
Note that de Cataldo and Migliorini proved the equality $(\ref{eq;7.2.1})$
holds as Chow motives {\rm (\cite{dcm})}.
\end{rem}

We prove the following equality.
\begin{thm}[Theorem \ref{thm;12.3.10}] \label{thm;12.2.6}
Let $X$ be a smooth surface over an algebraically closed field $k$,
and $D$ be a smooth divisor of  $X$.
Let $X^{[\vecv]}$ denote the parabolic Hilbert scheme
of points on $(X,D)$ with parabolic length $\vecv$ $(\vecv\in\nbiga)$.
In $K_0(V_k)\otimes\nbigr$,
the following formula holds:
\begin{equation} \label{eq;12.2.3}
\sum_{\vecv\in\nbiga}
 X^{[\vecv]}\cdot x^{\vecv}=
 \prod_{i=1}^{\infty} S_{i-1}(X,x_0^i)
\times
\prod_{\alpha<0}\prod_{i=1}^{\infty} S_{i-1}(D,x_0^i\!\cdot\! x_{\alpha})
\times
 \prod_{\alpha>0}
 \prod_{i=0}^{\infty}
 S_{i}(D,x_0^i\!\cdot\! x_{\alpha}).
\end{equation}
When the ground field $k$ is the complex number field $\cnum$,
we obtain the formula of the generating function
of the $E$-polynomials of  
$X^{[\vecv]}$
in terms of those of the symmetric powers of $X$ and $D$.
\mbox{{}}\hfill\qed
\end{thm}

As a corollary,
we can calculate the Betti number of
the parabolic Hilbert scheme
for a smooth projective surface $X$ and a smooth divisor $D$.
In general, we denote the Poincar\'e polynomial by $P(Y,z)$
with a variable $z$,
and we denote the $i$-th Betti number by $b_i(Y)$
for a topological space $Y$.
We denote the generating function
$\sum_{n=0}^{\infty} P(X^{[n]},z)\cdot t^n$
by $\tilde{P}_X(z,t)$, which was calculated by G\"ottsche:
\[
 \tilde{P}_X(z,t)=
 \prod_{m\geq 1}
 \frac{(1+z^{2m-1}t^m)^{b_1(X)}(1+z^{2m+1}t^m)^{b_3(X)}}
 {(1-z^{2m-2}t^m)^{b_0(X)}(1-z^{2m}t^m)^{b_2(X)}(1-z^{2m+2}t^m)^{b_4(X)}}.
\]
We can derive the following formula of the generating function
of the Poincar\'e polynomials of parabolic Hilbert schemes
of points
on $(X,D)$ from the equality (\ref{eq;12.2.3}).
For $a=0,1$ and for $\alpha\in\seisuu-\{0\}$,
we put as follows:
\[
 \tilde{P}_D(z,x_0,x_{\alpha},a):=
 \prod_{m\geq a}
 \frac{(1+z^{2(m-a)+1}x_0^mx_{\alpha})^{b_1(D)}}
       {(1-z^{2(m-a)}x_0^mx_{\alpha})^{b_0(D)}
        (1-z^{2(m-a)+2}x_0^mx_{\alpha})^{b_2(D)}}.
\]
\begin{cor}
\label{cor;7.6.1}
Let $X$ be a smooth projective surface over 
the complex number field $\cnum$,
and $D$ be a smooth divisor of $X$.
We consider the generating function of the Poincar\'e polynomials
with variable $z$
of the parabolic Hilbert schemes.
The following formula holds in $\seisuu[z]\otimes\nbigr$:
\begin{equation}  \label{eq;12.3.1}
\sum_{\vecv\in \nbiga}
 P(X^{[\vecv]},z)\cdot 
 x^{\vecv}
 =
 \tilde{P}_X(z,x_0)\times 
\prod_{\alpha<0}
 \tilde{P}_D(z,x_0,x_{\alpha},1)
\times 
\prod_{\alpha>0}
 \tilde{P}_D(z,x_0,x_{\alpha},0).
\end{equation}
\mbox{{}}\hfill\qed
\end{cor}

\begin{rem}
By using the result of de Cataldo-Migliorini {\rm \cite{dcm2}}
and our cell decomposition in section {\rm \ref{section;7.1.5}},
we can show that the equality $(\ref{eq;12.2.3})$
holds as the Chow motive.
The purpose of the argument
in subsection {\rm \ref{subsection;7.1.10}}
and section {\rm \ref{section;7.1.11}}
is to clarify the combinatorics
to derive the formula $(\ref{eq;12.3.1})$ directly
from the cell decomposition.
\end{rem}

\subsubsection{The incidence varieties and the induced operations}
\label{subsubsection;3.3.3}

For any element
$\vecu$ of $\tilde{C}$,
the variety $Y_{\vecu}$ is defined to be
$X\,\,(\rho(\vecu)= 0)$
or $D\,\,(\rho(\vecu)\neq 0)$.
For any element $\vecv\in \nbiga$,
we put
$\tilde{C}(\vecv):=\{\vecu\in\tilde{C}\,|\,\vecv+\vecu\in \nbiga\}$.
Then we define
the incidence variety
$Z(\vecv,\vecu)\subset
X^{[\vecv]}\times Y_{\vecu}\times X^{[\vecv+\vecu]}$
for any 
$\vecv\in \nbiga$ and $\vecu\in\tilde{C}(\vecv)$.
\label{subsub;7.8.10}
It is defined to be the closure of the following locally closed
subscheme:
\[
 \Bigl\{
  \bigl(
   (I_1,\nbigg_{1\,\ast}),x,(I_2,\nbigg_{2\,\ast})
  \bigr)
 \,\Bigl|\Bigr.\,
\begin{array}{l}
   (I_1,\nbigg_{1,\ast})=(I_2,\nbigg_{2,\ast})
 \mbox{ on $X-\{x\}$, }\\
 \length_x(\nbigo/I_1,\nbigg_{1,\ast})=\max(0,-\vecu),\,
 \length_x(\nbigo/I_2,\nbigg_{2,\ast})=\max(0,\vecu)
\end{array}
 \Bigr\}.
\]
The dimension of $Z(\vecv,\vecu)$
is $d(\vecv)+\rho_0(\vecu)+1+g(\vecu)$.

Let $H^{\ast}$ be an appropriate cohomology theory
defined for algebraic varieties.
For example, we consider the singular cohomology
defined for algebraic varieties over the complex number field $\cnum$.
We put as follows:
\[
 \hyperh(X,D):=\bigoplus_{\vecv\in \nbiga}
 H^{\ast}(X^{[\vecv]})\otimes x^{\vecv}
\]
Here $x^{\vecv}$ is added to emphasize the grading.
We give the multi-grading of $\hyperh(X,D)$
over $\seisuu_{\geq\,0}\times \nbiga$.
Let $a\otimes x^{\vecv}$ be an element
of $H^{\ast}(X^{[\vecv]})\otimes x^{\vecv}$.
The multi-degree $\deg(a\otimes x^{\vecv})$
is defined to be
$(\deg(a),\vecv)\in\seisuu_{\geq\,0}\times \nbiga$,
where $\deg(a)\in\seisuu_{\geq\,0}$ denotes the degree of $a$ in
$H^{\ast}(X^{[\vecv]})$.

From the incidence varieties $Z(\vecv,\vecu)$,
we obtain the operators
$\gminiq_{\vecu}(a):H^{\ast}(X^{[\vecv]})
     \lrarr H^{\ast}(X^{[\vecv+\vecu]})$:
We denote the projection of
$X^{[\vecv]}\times Y_{\vecu}\times X^{[\vecv+\vecu]}$
to the $i$-th component by $p_i$.
Then for any $c\in H^{\ast}(X^{[\vecv]})$,
the element $\gminiq_{\vecu}(a)(c)$
is defined by the correspondence as usual:
\[
 \gminiq_{\vecu}(a)(c):=
 p_{3\,\ast}\bigl(p_1^{\ast}(c)\cup p_2^{\ast}(a)\cup[Z(\vecv,\vecu)]\bigr).
\]
When $\vecu$ is not contained in $\tilde{C}(\vecv)$,
then $\gminiq_{\vecu}(a)(c)$ 
is defined to be $0$
for any element $c\in H^{\ast}(X^{[\vecv]})$.
As a result,
we obtain the operators
$\gminiq_{\vecu}(a):\hyperh(X,D)\lrarr\hyperh(X,D)$.
We have the grading of $End(\hyperh(X,D))$
over $\seisuu\times \TT$,
which is naturally induced by the grading of
$\hyperh(X,D)$.
Note that the multi-degree of the element
$\gminiq_{\vecu}(a)$ of $ End(\hyperh(X,D))$
is as follows:
\[
 \deg(\gminiq_{\vecu}(a))=
 \Bigl(\deg(a)+2(\rho_0(u)-1-g(-\vecu)),\vecu\Bigr)
 \in \seisuu\times\TT
\]

For any $\vecv\in \TT$,
the number $\epsilon_{\vecv}$ is defined to be $0\,\,(\vecv\neq 0)$
or $1\,\,(\vecv=0)$.
For each element $\vecu\in \tilde{C}$,
the number $sgn(\vecu)$ is defined to be
$1$ $(\vecu\in C)$
or $-1$ $(\vecu\in -C)$.
The next theorem is the main result of the paper.
\begin{thm}[Heisenberg relation, Corollary \ref{cor;3.10.5}]
Let $\vecu_1,\vecu_2$ be elements of $\tilde{C}$.
Let $a_i$ be an element of $H^{\ast}(Y_{\vecu_i})$
for $i=1,2$.
Then the following relation holds:
\[
 \bigl[\gminiq_{\vecu_1}(a_1),\gminiq_{\vecu_2}(a_2)\bigr]
=
 \epsilon_{\vecu_1+\vecu_2}\cdot
 \mu(\vecu_1)
 \cdot
 \left(
 \int_{Y_{\vecu_1}}a_1\cdot a_2\right)
 \cdot
 id_{\hyperh(X,D)}.
\]
Here $id_{\hyperh(X,D)}$ denotes
the identity of $\hyperh(X,D)$,
and we put as follows:
\[
 \mu(\vecu_1)=
 sgn(\vecu_1)\cdot
 \Bigl(-|\rho_0(\vecu_1)|+|\vecu_1|_- \Bigr)
  ^{\epsilon_{\rho(\vecu_1)}}
 (-1)^{\rho_0(\vecu_1)-|\vecu_1|_-}.
\]
In particular, the number $\mu(\vecu_1)$ is not $0$.
\hfill\qed
\end{thm}

\subsubsection{The representation of the Heisenberg algebra}

By the main theorem,
we obtain the structure of the representation
of the Heisenberg algebra on $\hyperh(X,D)$
as in the case of Hilbert schemes.
We follow the argument in \cite{l}.
We consider the following vector spaces:
\[
 W_+:=\bigoplus_{\vecu\in C} H^{\ast}(Y_{\vecu})\otimes x^{\vecu},
\quad
 W_-:=\bigoplus_{\vecu\in -C} H^{\ast}(Y_{\vecu})\otimes x^{\vecu},
\quad
 W:=W_+\oplus W_-.
\]
The multi-degree of  $W$ over the set
$\seisuu_{\geq\,0}\times \tilde{C}\subset
 \seisuu_{\geq\,0}\times \TT$
is defined by $\deg(a\otimes x^{\vecu})=\deg(\gminiq_{\vecu}(a))$.
We have the following non-degenerate graded skew-symmetric pairing
$\langle\cdot\,,\cdot\rangle$
over  $W$:
\[
 \bigl\langle\,a_1\otimes t^{\vecu_1},a_2\otimes t^{\vecu_2}\,\bigr\rangle=
 \epsilon_{\vecu_1+\vecu_2}\cdot
\mu(\vecu_1)\cdot
 \left(
 \int_{Y_{\vecu_1}}a_1\cdot a_2
 \right).
\]
We have the natural grading of
the tensor algebra $TW:=\bigoplus_n W^{\otimes\,n}$
over $\seisuu_{\geq\,0}\times \TT$
induced by the grading of $W$.
The two sided ideal $J$ of
$TW$
is generated by the expressions
$[v,w]-\langle v,w\rangle$ for any $v,w\in W$.
Then the oscillator algebra $\nbigh$
is defined to be $\nbigh:=TW/J$.
The correspondence
 $x\otimes t^{\vecu}\longmapsto
\gminiq_{\vecu}(x)$
induces the homomorphism
of $TW$ to $End\bigl(\hyperh(X,D)\bigr)$
of the graded associative algebras.
The theorem implies that it induces the morphism of
$\nbigh$
to $End\bigl(\hyperh(X,D)\bigr)$,
i.e.,
$\hyperh(X,D)$ is the $\nbigh$-representation.

Let $\nbigh\!\cdot\! W_-$ denote the left ideal generated by $W_-$.
Then we have the left $\nbigh$-module
$\nbigh/\nbigh\!\cdot\! W_-$,
which is naturally isomorphic to the
graded symmetric powers $\Sym(W_+)$ of
the space $W_+$.
Then we have the $\nbigh$-homomorphism
$\Psi:\nbigh/\nbigh\!\cdot\! W_-\lrarr \hyperh(X,D)$
induced by the correspondence $1\longmapsto 1_{H^0(X^{[0]})}$.
As usual, the module $\nbigh/\nbigh\!\cdot\! W_-$ is irreducible,
and the generating function of the dimensions
of $\nbigh/\nbigh\!\cdot \!W_-$ is same as
that for $\hyperh(X,D)$,
which is obtained in Corollary \ref{cor;7.6.1}.
Thus the morphism $\Psi$ is an isomorphism
as in the case of Hilbert schemes.
In all, we obtain the following theorem.
\begin{thm}
We have the $\nbigh$-representation structure
of $\hyperh(X,D)$.
It is isomorphic to
$\nbigh/\nbigh\cdot W_-$,
which is irreducible.
\hfill\qed
\end{thm}


\vspace{.1in}
\noindent
{\bf Acknowledgement}

The author expresses his gratitude
to Mark Andrea de Cataldo for the stimulating discussion.
He is also sincerely grateful to Luca Migliorini for his interest
to the parabolic Hilbert schemes.
The author is grateful to the colleagues in Osaka City University.
In particular, he thanks Mikiya Masuda for his encouragement.

The author thanks the financial supports by Japan Society for
the  Promotion of Science and the Sumitomo Foundation.
The paper was written during his stay at the Institute
for Advanced Study. The author is grateful to their excellent
hospitality. 
He also acknowledges National Scientific Foundation
for a grant DMS 9729992.


\section{Preliminary}

\subsection{Smoothness of parabolic Hilbert scheme}

Let $X$ be a smooth surface over an algebraically closed field $k$,
and $D$ be a smooth divisor of $X$.
Let $(I,\nbigg_{\ast})$ be a parabolic ideal
of points on $(X,D)$ 
with parabolic length $\vecv$ $(\vecv\in\nbiga)$.
We have the numbers $\alpha_-$ and $\alpha_+$
such that
$\rho_{\alpha}(\vecv)=0$ unless $\alpha_-\leq \alpha<\alpha_+$.
Let take a projective resolution of $I$,
that is,
$V_{\cdot}:=(V_{-1}\rarr V_0)\simeq I$.
We denote $V_{i}\otimes\nbigo_D$ by $V_{D,i}$.
Then we have the complexes
$C_1:=Cone(\nhom_X(V_{\cdot},V_{\cdot})\lrarr\nbigo)$ on $X$,
and $C_2:=\nhom_D(V_{D\,\cdot},V_{D\,\cdot})[1]$
on $D$.

We put $V_D^{(\alpha)}:=\ker(V_{D\,0}\lrarr \nbigg_{\alpha})$.
Note $V_D^{(\alpha)}=V_{D\,-1}$ for any $\alpha\leq \alpha_-$,
and $V_D^{(\alpha)}=V_{D\,0}$  for any $\alpha\geq \alpha_+$.
We put $A_{i,i+1}:=Cone(V_D^{(\alpha)}\rarr V_D^{(\alpha+1)})$.
Then we obtain the complexes
$\nhom_D(A_{\alpha,\alpha+1},A_{\alpha,\alpha+1})$
for any $\alpha$.
Thus we have the morphisms
$\nhom_D(V_D^{(\alpha)},V_D^{(\alpha)})
 \lrarr \nhom_D(V_D^{(\alpha)},V_D^{(\alpha+1)})$
and
$\nhom_D(V_D^{(\alpha+1)},V_D^{(\alpha+1)})
 \lrarr \nhom_D(V_D^{(\alpha)},V_D^{(\alpha+1)})$.
On the other hand,
we have the naturally defined morphisms
$\nhom_D(V_D^{(\alpha_+)},V_D^{(\alpha_-)})
 \lrarr \nhom_D(V_D^{(\alpha)},V_D^{(\alpha)})$
 for any $\alpha_-\leq\alpha\leq\alpha_+$.
Then we put as follows:
\[
 C_2:=
 \Bigl(
 \nhom_D(V_D^{(\alpha_+)},V_D^{(\alpha_-)})
 \lrarr
 \bigoplus_{i=\alpha_-}^{\alpha_+} 
 \nhom_D(V_{D}^{(\alpha)},V_{D}^{(\alpha)})
 \lrarr
 \bigoplus_{\alpha=\alpha_-}^{\alpha_+-1}
 \nhom_D(V_{D}^{(\alpha)},V_{D}^{(\alpha+1)})\Bigr).
\]
Here the degree of elements of
$\nhom_D(V_D^{(\alpha_+)},V_D^{(\alpha_-)})$ is $-2$.

We have the composite of the morphisms
$C_1\lrarr \nhom_X(V_{\cdot},V_{\cdot})[1]\lrarr C_3$.
On the other hand, 
the morphism $C_2\lrarr C_3$ is induced by
the naturally defined morphisms
$\nhom_D(V_{D}^{(\alpha)},V_D^{(\alpha+1)})
 \lrarr \nhom_D(V_{D\,-1},V_{D\,0})$
and the following projection:
\[
 \bigoplus_{i=\alpha_-}^{\alpha_+}
 \nhom_D(V_D^{(\alpha)},V_D^{(\alpha)})
 \lrarr \nhom_D(V_{D\,-1},V_{D\,-1})\oplus\nhom_D(V_{D\,0},V_{D\,0}).
\]
Thus we obtain the following complexes:
\[
 C(I,\nbigg_{\ast}):=
 Cone(C_1\oplus C_2\lrarr C_3)[-1],
\quad
 C_{par}(I,\nbigg_{\ast}):=
 Cone(C_2\lrarr C_3)[-1].
\]

\begin{lem}\label{lem;6.29.1}
The hypercohomology group $\hyperh^i(C(I,\nbigg_{\ast}))$ vanish
unless $i\neq 0$.
\end{lem}
\pf
We consider the cohomology sheaves $\nbigh^i(C(I,\nbigg_{\ast}))$.
Let $P$ be a point of $X$
such that $P$ is not contained in the supports
of the coherent sheaves $\nbigo/I$ or $\nbigk_{\alpha}$ $(\alpha\neq 0)$.
Then the complexes $C_i$ are acyclic around the point $P$.
Thus the support of $\nbigh^i(C(I,\nbigg_{\ast}))$ is $0$-dimensional.
Since the morphism
$\nhom_X(V_{-1},V_0)\lrarr \nhom_D(V_{D\,-1},V_{D\,0})$
is surjective,
we obtain $\nbigh^1(C(I,\nbigg_{\ast}))=0$.
It is easy to see that $\nbigh^i(C(I,\nbigg_{\ast}))=0$ for $i<0$.
Thus we obtain the result
by using the spectral sequence.
\hfill\qed


\begin{cor}\label{cor;smooth}
The parabolic Hilbert scheme points  on a smooth projective
surface is smooth.
\end{cor}
\pf
The obstruction theory of the moduli of oriented parabolic sheaves are given
in section 4 of \cite{m2}.
Note that a parabolic ideal is equivalent to a parabolic sheaf
of rank 1 whose determinant bundle is $\nbigo_X$.
Lemma \ref{lem;6.29.1} implies that
the $(-1)$-th cohomology sheaf of the obstruction complex vanishes.
Thus we obtain the result.
\hfill\qed

The following corollary is also obtained from the obstruction theory.
\begin{cor}
The tangent space at the points corresponding 
to a parabolic ideal $(I,\nbigg_{\ast})$
is naturally isomorphic to $\hyperh^0(C(I,\nbigg_{\ast}))$.
\hfill\qed
\end{cor}


\subsection{Some isomorphisms}

\label{subsection;7.4.10}

Let $\alpha_-$ and $\alpha_+$ be integers such that
$\alpha_-<0<\alpha_+$.
Let $\nbiga(\alpha_-,\alpha_+)$ denote
the set of elements $\vecv$ of $\nbiga$
such that
$\rho_{\alpha}(\vecv)=0$
unless $\alpha_-\leq\alpha<\alpha_+$.
Let take an integer $\beta\neq 0$
such that $\alpha_-\leq\beta<\alpha_+$.
We associate $\vecv'(\beta)$
to $\vecv\in \nbiga(\alpha_-,\alpha_+)$ as follows:
When $\beta>0$,
we put as follows:
\[
 \rho_{\alpha}(\vecv'(\beta)):=
\left\{
\begin{array}{ll}
 0, & (\beta\leq \alpha \mbox{ or } \alpha< \alpha_- -\alpha_+ +\beta),\\
 \rho_{\alpha}(\vecv), & (0<\alpha<\beta\mbox{ or } \alpha_-\leq\alpha<0),\\
 \rho_{0}(\vecv)+\sum_{\alpha'\geq \beta}\rho_{\alpha'}(\vecv), &(\alpha=0),\\
 \rho_{\alpha-\alpha_- +\alpha_+}(\vecv),
 & (\alpha_- -\alpha_+ +\beta\leq \alpha<\alpha_-).
\end{array}
\right.
\]
Then $\vecv'(\beta)$ is the element of
$\nbiga(\alpha_--\alpha_++\beta,\beta)$.
When $\beta<0$,
we put as follows:
\[
 \rho_{\alpha}(\vecv'(\beta)):=
\left\{
\begin{array}{ll}
 0, & (\alpha\leq \beta \mbox{ or } \beta+\alpha_+-\alpha_-\leq\alpha),\\
 \rho_{\alpha}(\vecv), &
 (\beta<\alpha<0\mbox{ or }0<\alpha<\alpha_+), \\
 \rho_0(\vecv)-\sum_{\alpha'\leq\beta}\rho_{\alpha'}(\vecv),
 & (\alpha=0),\\
 \rho_{\alpha-\alpha_+ +\alpha_- }(\vecv),
 & (\alpha_+\leq \alpha< \beta+\alpha_+-\alpha_-).
 \end{array}
\right.
\]
Then  $\vecv'(\beta)$ is the element of 
$\nbiga(\beta+1,\beta+\alpha_+-\alpha_-)$.

Let $(I,\nbigg_{\ast})$ be a parabolic ideal
with parabolic length $\vecv\in \nbiga(\alpha_-,\alpha_+)$.
We associate 
the parabolic ideal $(I',\nbigg_{\ast})$
with parabolic length $\vecv'(\beta)$ above
as follows:

When $\beta>0$,
we put $I':=\ker(I\lrarr\nbigg_{\beta})$.
The sheaf $I'$ is naturally an ideal sheaf of $0$-schemes
of length
$\rho_0(\vecv)+\sum_{\alpha'\geq \beta}\rho_{\alpha'}(\vecv)$.
It is easy to see that
we have the following exact sequence
of $\nbigo_D$-modules:
\[
 0\lrarr \nbigg_{\beta}\otimes\nbigo(-D)
 \lrarr I'\otimes\nbigo_D\lrarr \nbigl_{\beta}\lrarr 0.
\]
Here $\nbigl_{\beta}$ denotes the kernel
of the morphism $I_D\lrarr\nbigg_{\beta}$.
We have the naturally induced filtrations
on the quotient $\nbigg_{\beta}$
and the subsheaf $\nbigl_{\beta}$ of $I\otimes\nbigo_D$.
By taking appropriate shift of the indices of the filtration
of $\nbigg_{\beta}$,
we obtain the desired filtration on $I'\otimes\nbigo_D$.

When $\beta<0$,
we put $I':=\ker(I\lrarr\nbigg_{\beta})\otimes\nbigo(D)$.
The sheaf $I$ is naturally contained in $I'$,
and the cokernel is isomorphic to 
$\nbigl_{\beta}\otimes\nbigo(D)=
 \ker\bigl(I\otimes\nbigo_D\lrarr \nbigg_{\beta}\bigr)
  \otimes \nbigo(D)$.
It is easy to see that we have the following exact sequence:
\[
 0\lrarr \nbigg_{\beta}\lrarr I'\otimes\nbigo_D
 \lrarr \nbigl_{\beta}\otimes\nbigo_D\lrarr 0.
\]
We have the induced filtrations
on $\nbigg_{\beta}$ and $\nbigl_{\beta}\otimes\nbigo_D$.
By taking the appropriate shift of the indices
of the filtration of $\nbigl_{\beta}$,
we obtain the desired filtration on $I'\otimes\nbigo_D$.

The correspondence $(I,\nbigg_{\ast})$
to $(I',\nbigg'_{\ast})$ above
induces the isomorphism $\Psi_{\vecv,\vecv'(\beta)}$ of
$X^{[\vecv]}$
to $X^{[\vecv'(\beta)]}$.


\begin{rem} \label{rem;7.5.2}
In particular,
if $\beta=-1$ then $\vecv'(\beta)$ is contained
in $\TT_0\oplus\TT_+$.
\hfill\qed
\end{rem}

\label{subsection;7.5.1}


\subsection{Symmetric powers}

\label{subsection;7.1.10}
\subsubsection{$\nbigs(\Lambda)$ and $\nbigt(\Lambda)$}

For any set $\Lambda$,
we put as follows:
\begin{equation} \label{eq;nbigs}
\begin{array}{l}
 \nbigs(\Lambda):=
 \Bigl\{f:\Lambda\lrarr \seisuu_{\geq\,0},\,
 \#\{\lambda\in\Lambda\,|\,f(\lambda)\neq 0\}<\infty
 \Bigr\}.\\
 \mbox{}\\
 \nbigt(\Lambda):=
 \Bigl\{f:\Lambda\lrarr \seisuu_{\geq\,0},\,
 0<\#\{\lambda\in\Lambda\,|\,f(\lambda)\neq 0\}<\infty
 \Bigr\}.\\
\end{array}
\end{equation}
We have the natural injection
$\nbigt(\Lambda)\lrarr \nbigs(\Lambda)$.

For any element $\rho\in\nbigs(\Lambda)$,
we put as follows:
\[
 |\rho|:=\sum_{\lambda\in \Lambda}\rho(\lambda)\in\seisuu_{\geq \,0}
\]
Thus we obtain the homomorphism
$|\cdot|:\nbigs(\Lambda)\lrarr \seisuu_{\geq\,0}$
of semi-groups.
We denote the restriction of the map $|\,\cdot\,|$
to $\nbigt(\Lambda)$ by the same notation.

For any element $f\in \nbigs(\nbigs(\Lambda))$,
we put as follows:
\[
 \phi_{\Lambda}(f)
 :=\sum_{\rho\in\nbigs(\Lambda)}f(\rho)\rho \in \nbigs(\Lambda)
\]
Thus we obtain
the homomorphism
$\phi_{\Lambda}:\nbigs(\nbigs(\Lambda))\lrarr\nbigs(\Lambda)$
of semi-groups.
We denote the restriction of the morphism $\phi_{\Lambda}$
to $\nbigs(\nbigt(\Lambda))$ 
by the same notation.
Let $\pi$ be an element of $\nbigs(\Lambda)$.
Then the subset $\phi_{\Lambda}^{-1}(\pi)$ of $\nbigs(\nbigt(\Lambda))$
is denoted by $\nbigs(\nbigt(\Lambda),\pi)$.

When $\Lambda$ consists only one element,
then $\nbigs(\Lambda)$ is isomorphic to $\seisuu_{\geq\,0}$
and $\nbigt(\Lambda)$ is isomorphic to $\seisuu_{\geq \,1}$.
We will use the identification without mentioning.
In this case,
$\phi_{\Lambda}(f)$ is same as the number
$\sum_{a\in\seisuu_{\geq\,1}}f(a)\cdot a\in
 \nbigs(\Lambda)\simeq\seisuu_{\geq\,0}$
for any $f\in \nbigs(\nbigt(\Lambda))$.

\subsubsection{Stratification of symmetric power}

Let $Y$ be a variety over $k$.
We denote the $n$-th symmetric power
of $Y$ by $Y^{(n)}$.
We have the projection $\pi_1:Y^n\lrarr Y^{(n)}$.
We put $Y^n_{\ast}:=\{(x_1,\ldots,x_n)\,|\,x_i\neq x_j\,(i\neq j)\}$
and $Y^{(n)}_{\ast}:=\pi_1(Y^n_{\ast})$.

Let $\Lambda$ be a set.
For an element $\rho\in \nbigs(\Lambda)$,
we put $Y^{(\rho)}:=\prod_{\lambda\in\Lambda} Y^{(\rho(\lambda))}$.
We have the naturally defined morphism
$\pi_2:Y^{(\rho)}\lrarr Y^{(|\rho|)}$.
We put $Y^{(\rho)}_{\ast}:=\pi_2^{-1}(Y^{|\rho|}_{\ast})$.
Assume that variety
$Z_{\lambda}$ is  given for each $\lambda\in \Lambda$.
Then we have the canonical projection
$\pi_3:\prod_{\lambda\in\Lambda}(Y\times Z_{\lambda})^{(\rho(\lambda))}
 \lrarr Y^{(\rho)}$.
We put as follows:
\[ 
 \Bigl(
 \prod_{\lambda\in\Lambda}
 (Y\times Z_{\lambda})^{(\rho(\lambda))}
 \Bigr)_{\ast}
 :=
 \pi_3^{-1}(Y^{(\rho)}_{\ast}).
\]

We put $N=\nbigs(\Lambda)$ and $N'=\nbigt(\Lambda)$.
Let $\rho$ be an element of $N$
and $f$ be an element of $\nbigs(N',\rho)$.
For any element $a\in N'$ and $\lambda\in \Lambda$,
we have the diagonal morphism $Y\lrarr Y^{a(\lambda)}$.
It induces the morphism
$Y^{f(a)}\lrarr (Y^{a(\lambda)})^{f(a)}=Y^{a(\lambda)\cdot f(a)}$.
Thus we obtain the morphism
$Y^{f(a)}\lrarr \prod_{\lambda\in \Lambda}Y^{a(\lambda)\cdot f(a)}$.
Then we obtain the following morphism:
\[
 \prod_{a\in N'}Y^{f(a)}\lrarr
 \prod_{a\in N'}
 \prod_{\lambda\in\Lambda} Y^{a(\lambda)\cdot f(a)}
\lrarr
 \prod_{\lambda\in\Lambda}Y^{\phi(f)(\lambda)}
=\prod_{\lambda\in\Lambda}Y^{\rho(\lambda)}.
\]
Thus we obtain the morphism
$Y^{(f)}\lrarr Y^{(\rho)}$.
The following lemma is easy to see.
\begin{lem}
The following morphism gives the decomposition
of $Y^{(\rho)}$:
\[
 \coprod_{f\in\nbigs(N',\rho)} Y_{\ast}^{(f)}
 \lrarr 
 Y^{(\rho)}.
\]
In particular,
it gives the equality in the Grothendieck ring $K_0(V_k)$
of smooth $k$-varieties.
\hfill\qed
\end{lem}

\subsubsection{Some formulas in $K_0(V_k)[[N]]$}

Let $\Lambda$ be a finite set, and $N$ be $\nbigs(\Lambda)$.
We denote the semi-group ring associated to $N$ by
$\seisuu[N]$,
i.e.,
$\seisuu[N]$ is the ring generated
by the elements $\{x^{\rho}\,|\,\rho\in N\}$
with the relations
$x^{\rho_1+\rho_2}=x^{\rho_1}\cdot x^{\rho_2}$
$(\rho_1,\rho_2\in N)$ and $x^0=1$.
We have the ideal $\gminip$ generated
by $\{x^{\rho}\,|\,\rho\in N-\{0\}\,\,\}$.
We denote the completion with respect to $\gminip$-adic topology
by $\seisuu[[N]]$.
Let $\delta_{\lambda}$ denote the element of $N$
such that $\delta_{\lambda}(\lambda')$ is $1$ $(\lambda=\lambda')$
or $0$ $(\lambda\neq \lambda')$.
We put $x_{\lambda}:=x^{\delta_{\lambda}}$.
The ring $\seisuu[N]$ is isomorphic to a polynomial ring
generated by $\{x_{\lambda}\,|\,\lambda\in\Lambda\}$
over $\seisuu$.
The ring $\seisuu[[N]]$ is isomorphic to the ring
of formal power series generated by 
$\{x_{\lambda}\,|\,\lambda\in\Lambda\}$
over $\seisuu$.

Let $\Lambda$ be a countable set.
For any finite subset $\Lambda_i$ of $\Lambda$,
we put $N_i=\nbigs(\Lambda_i)$.
For $\Lambda_1\subset\Lambda_2$,
we have the natural morphism
$N_2\lrarr N_1$ and thus $\seisuu[[N_2]]\lrarr \seisuu[[N_1]]$.
We put $\seisuu[[N]]:=\varprojlim_{\Lambda_i\subset\Lambda} \seisuu[[N_i]]$.

We denote the Grothendieck ring of the varieties 
over $k$ by $K_0(V_k)$.
Then we denote the ring $K_0(V_k)\otimes_{\seisuu} \seisuu[[N]]$
by $K_0(V_k)[[N]]$.
Let $Z_1,\ldots Z_n$ and $Y$ be varieties.
The following lemma is easy to see.
\begin{lem}
The following equality holds in $K_0(V_k)[[N]]$:
\[
 \prod_{i=1}^n
 \Bigl(\sum_{\rho\in N}
  (Y\times X_i)^{(\rho)}x^{\rho}
 \Bigr)_{\ast}
=
 \sum_{\rho\in N}
 \Bigl(Y\times
 \sum_{i=1}^n X_i
 \Bigr)^{(\rho)}_{\ast}.
\]
\hfill\qed
\end{lem}

Let $\Lambda_i\,(i=1,2)$ be sets.
We put $N_i:=\nbigs(\Lambda_i)$
and $N_i':=\nbigt(\Lambda_i)$.
Let $g:N_1\lrarr N_2$ be a homomorphism of semi-groups.
\begin{lem}
Let $Y$ be a variety.
Assume that a variety $Z_{\rho}$ is given for each $\rho\in N_1$.
Then we have the following equality:
\begin{equation} \label{eq;6.30.15}
 \sum_{f\in \nbigs(N_2')}
 \Bigl[
 \prod_{\mu\in N'_2}
 \Bigl(
 \sum_{\rho\in g^{-1}(\mu)}Y\times Z_{\rho}
 \Bigr)^{(f(\mu))}
 \Bigr]_{\ast}
 \cdot
 x^{\phi_{\Lambda_2}(f)}
=
 \sum_{F\in \nbigs(N_1')}
 \Bigl[
 \prod_{\rho\in N_1'}(Y\times Z_{\rho})^{(F(\rho))}
 \Bigr]_{\ast}\cdot x^{g(\phi_{\Lambda_1}(F))}.
\end{equation}
\end{lem}
\pf
The left hand side can be rewritten as follows:
\begin{multline}
 \Bigl[
  \prod_{\mu\in N_2'}
 \Bigl(
 \sum_{j=0}^{\infty}
 \bigl(\sum_{\rho\in g^{-1}(\mu)}Y\times Z_{\rho}
 \bigr)^{(j)} \cdot x^{j\cdot \mu}
 \Bigr)
 \Bigr]_{\ast}
=
 \Bigl[
 \prod_{\mu\in N_2}
 \prod_{\rho\in g^{-1}(\mu)}
 \Bigl(
 \sum_{j=0}^{\infty}
 (Y\times Z_{\rho})^{(j)}
 x^{j\cdot g(\rho)}
 \Bigr)
 \Bigr]_{\ast}\\
 =\Bigl[
 \prod_{\rho\in N_1}
 \Bigl(
 \sum_{j=0}^{\infty}
 (Y\times Z_{\rho})^{(j)}
 x^{j\cdot g(\rho)}
 \Bigr)
 \Bigr]_{\ast}
\end{multline}
It can be rewritten as follows:
\[
 \sum_{F\in \nbigs(N_1')}
 \Bigl[
 \prod_{\rho\in N_1'}
 (Y\times Z_{\rho})^{(F(\rho))}
 \Bigr]_{\ast}
 \cdot x^{g(\sum F(\rho)\rho)}.
\]
By definition,
we have the equality $\sum F(\rho)\rho=\phi_{\Lambda_1}(F)$.
Thus we are done.
\hfill\qed

For elements $\lambda\in \Lambda_1$
and $\lambda_2\in\Lambda_2$,
we have the numbers $a(\lambda,\lambda_2)\in\seisuu_{\geq\,0}$
satisfying the equality
$g(\delta_{\lambda})=
 \sum_{\lambda_2\in\Lambda_2}a(\lambda,\lambda_2)\delta_{\lambda_2}$.
Let $c:N_1\lrarr \seisuu_{\geq \,0}$ be a homomorphism
of semi-groups.
For an element $\lambda\in \Lambda_1$,
we put $c_{\lambda}:=c(\delta_{\lambda})$.
Let $A^l$ denote an $l$-dimensional affine space over $k$.

\begin{lem}\label{lem;6.30.20}
We have the following equality
in $K_0(V_k)[[N]]$:
\[
  \sum_{f\in \nbigs(N_2')}
 \Bigl[
 \prod_{\mu\in N'_2}
 \Bigl(
 \sum_{\rho\in g^{-1}(\mu)}Y\times A^{c(\rho)}
 \Bigr)^{(f(\mu))}
 \Bigr]_{\ast}
 \cdot
 x^{\phi_{\Lambda_2}(f)}
=
 \prod_{\lambda\in\Lambda_1}
 \Bigl(
\sum_{j=0}^{\infty}
 (Y \times A^{c_{\lambda}})^{(j)}\cdot z(\lambda)^j
 \Bigr).
\]
Here we put
$z(\lambda):=x^{g(\delta_{\lambda})}=\prod_{\lambda_2\in \Lambda_2}
 x_{\lambda_2}^{a(\lambda,\lambda_2)}$.
\end{lem}
\pf
We substitute $Z_{\rho}=A^{c(\rho)}$
into the equality (\ref{eq;6.30.15}).
Recall that
we have the equality $(Y\times A^a)^{(n)}=Y^{(n)}\times A^{an}$
in $K_0(V_k)$
by the argument due to Totaro (See \cite{g2}).
Thus the left hand side can be rewritten as follows:
\[
 \sum_{F\in\nbigs(N_1')}
 \Bigl[
 \prod_{\rho\in N_1'}
 Y^{(F(\rho))}
 \Bigr]_{\ast}
 \times A^{\sum F(\rho)\cdot c(\rho)}
 \cdot x^{g(\phi_{\Lambda_1}(F))}.
\]
It can be rewritten as follows:
\[
 \sum_{\eta\in N_1}
 \sum_{F\in \nbigs(N_1',\eta)}
 \Bigl[
 \prod_{\rho\in N_1'}
 Y^{(F(\rho))}
 \Bigr]_{\ast}
 \times A^{c(\eta)}
 \cdot x^{g(\eta)}
=\sum_{\eta\in N_1}
 Y^{(\eta)}\times A^{c(\eta)}\cdot x^{g(\eta)}
=
 \prod_{\lambda\in\Lambda_1}
 \Bigl(
\sum_{j=0}^{\infty}
 (Y \times A^{c_{\lambda}})^{(j)}\cdot z(\lambda)^j
 \Bigr).
\]
Thus we are done.
\hfill\qed

\subsubsection{$\nN$, $\mM$, $\nnN$ and $\mmM$}

\label{subsub;7.1.50}

We put $\nN:=\nbigs(\LLA)$ and
$\nN':=\nbigt(\LLA)$.
Note that we have the natural inclusions
$\nbiga\subset\nN\subset\TT$.
We put $\mM:=\nbigs(\LLB)$
and $\mM':=\nbigt(\LLB)$.
An element of $\mM$ or $\mM'$ is described as
$\eta=(\eta_{\alpha}\,|\,\alpha\in\LLA)$,
where $\eta_{\alpha}$ is an element of $\nbigs(C_{\alpha})$
(See subsubsection \ref{subsubsection;7.6.10}).
For any element $\eta=(\eta_{\alpha}\,|\,\alpha\in \LLA)$
of $\mM$,
we put $||\eta||:=\sum_{\alpha\leq 0}|\eta_{\alpha}|$.

By construction,
we have the maps $\phi_{\LLA}:\nbigs(\nN')\lrarr \nN$
and $\phi_{\LLB}:\nbigs(\mM')\lrarr \mM$.
We put $\nbigs(\nN',\vecv):=\phi_{\LLA}^{-1}(\vecv)$
for $\vecv\in \nN$.
We also put
$\nbigs(\mM',\eta):=\phi_{\LLB}^{-1}(\eta)$
for $\eta\in \mM$.
On the other hand, we have the map $\Psi:\mM\lrarr \nN$
defined as follows:
\begin{equation}\label{eq;Psi}
\Psi(\eta):=\sum_{\vecu\in C}\eta(\vecu)\cdot\vecu
 \in\nbiga\subset \nN.
\end{equation}
For an element $\vecv$ of $\nN$,
we put $\mM(\chi):=\Psi^{-1}(\vecv)$.
Note that the image of $\mM$ via the morphism $\Psi$
is $\nbiga$.

We put $\nnN:=\nbigs(\{0\})=\seisuu_{\geq\,0}$
and $\nnN':=\nbigt(\{0\}) =\seisuu_{\geq\,1}$.
We put $\mmM:=\nbigs(\LLB_0)$
and $\mmM':=\nbigt(\LLB_0)$.
By construction,
we have the maps $\phi_{\{0\}}:\nbigs(\nnN')\lrarr \nnN$
and $\phi_{\LLB_0}:\nbigs(\mmM')\lrarr \mmM$.
We put $\nbigs(\nnN',a):=\phi_{\{0\}}^{-1}(a)$
for $a\in \nnN$.
We also put
$\nbigs(\mmM',\rho):=\phi_{\LLB_0}^{-1}(\rho)$
for $\rho\in \mmM$.
Since $\mmM$ is isomorphic to $\nbigs(\nbigt(\{0\}))$,
we have the morphism $\phi_{\{0\}}:\mmM\lrarr \nnN$.
For an element $a\in \nnN$,
we put $\mmM(a):=\phi^{-1}_{\{0\}\,\ast}(a)$.


\section{ The cell decomposition of the punctual parabolic Hilbert
 schemes}

\label{section;7.1.5}

\subsection{Fixed point set}
We consider the case $X=A^2=\Spec k[x,y]$ and
$D=\{(x,y)\in X\,|\,x=0\}$.
We denote the point $(0,0)$ by $O$.
Let $\vecv$ be an element of $\nbiga$.
We consider the punctual parabolic Hilbert scheme
$X^{[\vecv]}_p$
which is the moduli of the parabolic ideals
$(I,\nbigg_{\ast})$
such that the union of the supports of the torsion sheaves
$\nbigo/I$ and $\nbigk_{\alpha}\,(\alpha\neq 0)$
is the set $\{O\}$.
In this section, we see the cell decomposition of 
the parabolic Hilbert schemes of points on the surface
by using the argument due to Ellingsrud and Str\o{}mme
\cite{es}.

We have the $2$-dimensional torus
$G_m^2=\Spec k[\lambda,\lambda^{-1},\mu,\mu^{-1}]$-action
on $X$ which preserves the divisor $D$:
$t\cdot(x,y)=(\lambda(t)\cdot x,\mu(t)\cdot y)$.
Thus we obtain the
$2$-dimensional torus action on 
$X^{[\vecv]}$ if $X=A^2$ and $D=A^1$.

The fixed point set 
of the $G_m^2$-action on the punctual parabolic Hilbert scheme
$X^{[\vecv]}_p$
is discrete and labeled by the elements
$\eta\in \mM(\vecv)$.
The correspondence is explained in the following:

Let $(I,\nbigg_{\ast})$ be the ideal
corresponding to a fixed point set.
The ideal $I$ is generated by the monomials
$x^jy^{b_j}$, $(j=0,1,\ldots,r)$
where $b_0\geq b_1\geq \cdots\geq b_r\geq 0$.
We put $b_i=0$ for any $i\geq r+1$.
The tuple $\{b_0,b_1,b_2,\ldots\}$ satisfies 
the equality $\sum b_i=n$.

We put $\nbige_{-1}:=\bigoplus_{j=1}^r \nbigo_X\cdot f_j$,
and $\nbige_{0}:=\bigoplus_{j=0}^r\nbigo_X\cdot v_j$.
The $G_m^2$-action on $\nbige_{-1}$ and $\nbige_{0}$
is given as follows:
$t\cdot v_j:=\lambda(t)^j\mu(t)^{b_{j}}v_j$ and
$t\cdot f_j:=\lambda(t)^j\mu(t)^{b_{j-1}}f_j$.
Then we obtain the $G_m^2$-equivariant resolution
$(\nbige_{-1}\rarr\nbige_0)$ of $I$.
Here the element $f_j$ is mapped to
$y^{b_{j-1}-b_j}v_j-x v_{j-1}\in \nbige_0$,
and the element $v_j$ is mapped to $x^jy^{b_j}\in I$.
We put $I_{D\,0}:=\nbigo_D\cdot v_0$
and $I_{D\,j}:=(\nbigo_D/y^{b_{j-1}-b_j}\nbigo_D)\cdot v_j$.
Then we have the natural isomorphism
$I_D\simeq
 \bigoplus_{j\geq 0} I_{D\,j}$.

Since $(I,\nbigg_{\ast})$ has the $G_m^2$-action,
the filtration $\nbigg_{\ast}$
is the direct sum of the filtrations
$\nbigg_{j\,\ast}=\{\nbigg_{j\,\alpha}\,|\,\alpha\in\seisuu\}$
of $I_{D\,j}\,(j=0,1,2,\ldots)$.
Here $\nbigg_{j\,\alpha}$ is a quotient of $\nbigg_{j\,\alpha-1}$,
we have $\nbigg_{j\,\alpha}=I_{D\,j}$ for any sufficiently small
$\alpha$,
and we have $\nbigg_{\alpha}=\bigoplus_j\nbigg_{j\,\alpha}$.
We put
$\nbigk_{j,\alpha}:=\ker(\nbigg_{j,\alpha}\lrarr\nbigg_{j,\alpha+1})$.
Note that
$\nbigk_{0,0}\simeq \nbigo_D$,
that the others
$\nbigk_{j,\alpha}$ $(j,\alpha)\neq(0,0)$ are torsion,
and that $\nbigk_{0,\alpha}=0$ for any
$\alpha<0$.

For a positive integer $\alpha>0$,
we put
$\eta_{\alpha}(j):=\length(\nbigk_{j,\alpha})$
for any $j\geq 0$.
For a non-positive integer $\alpha\leq 0$,
we put
$\eta_{\alpha}(j):=\length(\nbigk_{j,\alpha})$
for any $j>0$.
Then
$\eta_{\alpha}$ gives the element of  $\nbigs(C_{\alpha})$.
Here we used the identification
$C_{\alpha}$ with $\nbigs(\nbigs(\{\alpha\}))$
(resp. $\nbigs(\nbigt(\{\alpha\}))$)
for $\alpha>0$ (resp. $\alpha\geq 0$).

Thus we obtain the element
$\eta=(\eta_{\alpha}\,|\,\alpha\in\LLA)$ of $\mM$.
By definition, we have the equalities
$ 
 |\eta_{\alpha}|=\rho_{\alpha}(\vecv) $ for any 
$\alpha\in\LLA-\{0\}$.
We also have the equality
$\sum_{\alpha}\phi_{\{\alpha\}}(\eta_{\alpha})=
 \sum_{i=0}^{\infty}i\cdot (b_{i-1}-b_i)=
 \rho_0(\vecv)$.
Here we used the identification of
$C_{\alpha}$ with $\nbigs(\{\alpha\})$ $(\alpha>0)$,
or $\nbigt(\{\alpha\})$ $(\alpha\leq 0)$.
Thus it holds that $\Psi(\eta)=\vecv$.
Hence $\eta$ is the element of $\mM(\vecv)$.

On the other hand,
we can reconstruct the fixed point
from an element of $\mM(\vecv)$.
Thus we obtain the correspondence.
We denote the fixed point corresponding to $\eta\in \mM(\vecv)$
by $FP(\eta)$.

\subsection{The weight decomposition of the tangent space}

Let $\vecv$ be an element of $\nbiga$ 
and $\alpha_-$ be an integer.
Assume that $\rho_{\alpha}(\vecv)=0$ for any $\alpha<\alpha_-$.
Let $\eta$ be an element of $\mM(\vecv)$.
We have the $G_m^2$-action on the tangent space of
$X^{[\vecv]}$
at the fixed point $FP(\eta)$.
We calculate the weight decomposition of the $G_m^2$-representation.
Recall that the representation ring of $G_m^2$
is the Laurent polynomial ring
$\seisuu[\lambda,\lambda^{-1},\mu,\mu^{-1}]$.
It will be convenient to admit some infinite sums.
Thus we use the ring
$\seisuu[[\lambda,\mu]][\lambda^{-1},\mu^{-1}]$,
whose element is of the form
$\sum_{i,j\geq N}a_{i\,j}\cdot\lambda^i\cdot \mu^j$
for some $N\in\seisuu$ and $a_{i\,j}\in \seisuu$.
The ring is called the representation ring
for abbreviation.

For the given $\eta$,
we put $b_{i-1}:=\sum_{j\geq i}\sum_{\alpha} \eta_{\alpha}(j)$.
We also put
$\bar{h}_{\alpha,i}:=\sum_{\beta\geq\alpha}\eta_{\beta}(i)$
for $i>0$ or for $i=0,\,\alpha>0$.
Let $(I,\nbigg_{\ast})\in X^{[\vecv]}$
be a parabolic ideal corresponding
to the fixed point $FP(\eta)$.
Then the number $\bar{h}_{\alpha,i}$
is same as the length of the torsion sheaf $\nbigg_{\alpha,i}$.

\begin{prop} \label{prop;12.3.15}
In the representation ring of $G_m^2$,
the following equality holds:
\begin{multline} \label{eq;oct.22}
 T_{FP(\eta)}X^{[\vecv]}
 =
\sum _{1\leq i\leq j}\sum_{s=b_j}^{b_{j-1}-1}
 (\lambda^{i-j-1}\mu^{b_{i-1}-s-1}+\lambda^{j-i}\mu^{s-b_{i-1}})\\
+\sum_{
 (\alpha,i,j)\in S
  }
\sum_{a=0}^{\bar{h}_{\alpha,j}-1}
  \lambda^{j-i}
 \bigl(
   \mu^{b_j-b_i-\bar{h}_{\alpha,i}+a}
  -\mu^{b_j-b_i-\bar{h}_{\alpha-1,i}+a}
 \bigr) \\
-
\sum_{\alpha_-\leq\alpha\leq 0}\sum_{i\geq 1}
 \sum_{a=0}^{\eta_{\alpha-1}(i)-1}
 \lambda^{-i}\mu^{b_0-b_i+a-\bar{h}_{\alpha-1,i}} 
+
 \sum_{j\geq 0}
 \sum_{a=0}^{\bar{h}_{1,j}-1}
 \lambda^j
 \mu^{b_j-b_0-\bar{h}_{1,0}+a}.
\end{multline}
Here $S$ denotes the set
of $(\alpha,i,j)$
satisfying one of the following conditions:
$(i)$ $\alpha\geq \alpha_-$, $i\geq 1$, and $j\geq 1$,
$(ii)$ $\alpha>0$, $i\geq 1$ and $j=0$,
or
$(iii)$ $\alpha>1$, $i=0$ and $j\geq 0$.
\end{prop}
\pf
We have the natural isomorphism
$T_{FP(\eta)}X^{[\vecv]}
\simeq\hyperh^0(C(I,\nbigg_{\ast}))$.
We know that $\hyperh^i(C(I,\nbigg_{\ast}))$ vanishes unless $i=0$.
Thus we have the following $K$-theoretic equality:
\[
 \hyperh^0(C(I,\nbigg_{\ast}))=
 \hyperh^{\ast}(C_1)+\hyperh^{\ast}(C_{par})
 :=\sum (-1)^i\hyperh^{i}(C_1)+\sum (-1)^i\hyperh^i(C_{par}).
\]
Here we know that $\hyperh^{i}(C_1)$ vanishes
unless $i=0$,
and $\hyperh^0(C_1)$ is naturally isomorphic to
the tangent space of $X^{[n]}$
at the fixed point corresponding to $I$.
Thus the contribution of
$\hyperh^{\ast}(C_1)=\hyperh^0(C_1)$
was calculated in \cite{es}.
The result is the first term in the right hand side of the equality
(\ref{eq;oct.22}).

To see the contribution of $\hyperh^{\ast}(\nbigc_2)$,
we use the $G_m^2$-equivariant resolution $\nbige_{-1}\lrarr\nbige_0$
of $I$.
We have the direct sum decomposition
$\nbigg_{\alpha}=\bigoplus_{i}\nbigg_{\alpha,i}$ where $\nbigg_{\ast,i}$
is the filtration of the $I_{D\,i}$.
We put
 $E^{\alpha}
  :=\ker(\nbige_0\otimes\nbigo_D\lrarr \nbigg_{\alpha})$.
We have the natural decomposition
$E^{\alpha}=\bigoplus_i E^{\alpha}_i$.
Note that $E^{\alpha}=0$ for any sufficiently small $\alpha$
and $E^{\alpha}=\nbige_0\otimes\nbigo_D$
for any sufficiently large $\alpha$.

For $G_m^2$-equivariant sheaves $F_a$ $(a=1,2)$ over $D=A^1$,
let $Hom(F_1,F_2)$ denote the $G_m^2$-representation
$H^0(D,\nhom_{D}(F_1,F_2))$.
Although it is infinite dimensional representation,
it gives the element of our representation ring.
It can be checked by a formal calculation
that the contribution of $C_2$ is same as the following
in the representation ring of $G_m^2$:
\[
  \sum_{\alpha\geq \alpha_-}
   Hom(E^{\alpha},\nbigg_{\alpha})
-
 \sum_{\alpha\geq\alpha_-}
   Hom(E^{\alpha-1},\nbigg_{\alpha})=
\sum_{\alpha\geq 1}\sum_{j\geq 0}\sum_{i\geq 0}
\Bigl(
 Hom(  E^{\alpha}_{i}, \nbigg_{\alpha,j})
-
 Hom(E^{\alpha-1}_i, \nbigg_{\alpha,j})\Bigr).
\]
\begin{lem}
In the representation ring,
we have the following equality:
\[
 \sum_{i\geq 1}\sum_{j\geq 1}
 Hom(  
 E_{i}^{\alpha},
         \nbigg_{\beta,j})
=
  \sum_{i\geq 1}\sum_{j\geq 1}\sum_{a=0}^{\bar{h}_{\beta,j}-1}
 \lambda^{j-i}\mu^{b_j-b_i-\bar{h}_{\alpha,i}+a}.
\]
\end{lem}
\pf
The sheaf $E_i^{\alpha}$ is a locally free sheaf of rank 1
generated by the section with the weight
$\lambda^{i}\mu^{b_i+\bar{h}_{\alpha,i}}$.
The $\nbigg_{\beta,j}\,(j>0)$ is the torsion sheaf
of length $\bar{h}_{\beta,j}$ generated by
the section with the weight $\lambda^j\mu^{b_j}$.
Thus we obtain the result.
\hfill\qed

\vspace{.1in}

Thus the contribution of
$
 \sum_{\alpha\geq \alpha_-}\sum_{i\geq 1}\sum_{j\geq 1}
\bigl(
 Hom(
 E^{\alpha}_{i},
 \nbigg_{\alpha,j})
-
 Hom(
 E^{\alpha-1}_i, \nbigg_{\alpha,j})
\bigr)$
is as follows:
\[
 \sum_{\alpha\geq \alpha_-}\sum_{i\geq 1}\sum_{j\geq 1}
 \sum_{a=0}^{\bar{h}_{\alpha,j}-1}
\lambda^{j-i}
\bigl(\mu^{b_j-b_i-\bar{h}_{\alpha,i}+a}
-\mu^{b_j-b_i-\bar{h}_{\alpha-1,i}+a}
\bigr).
\]

Next we see the contribution of the $i=0$.
\begin{lem}
Assume that $\alpha\leq 0$.
In the representation ring,
we have the following equality:
\[
 \sum_{i\geq 1}
\bigl(
 Hom(
 E_{i}^{\alpha},\nbigg_{\alpha,0})
 -
 Hom(
 E_{i}^{\alpha-1},
     \nbigg_{\alpha,0})\bigr)=
 -\sum_{i\geq 1}\sum_{a=0}^{\eta_{\alpha-1}(i)-1}
 \lambda^{-i}\mu^{b_0-b_i+a-\bar{h}_{\alpha-1,i}}.
\]
\end{lem}
\pf
When $\alpha\leq 0$,
the sheaf $\nbigg_{\alpha,0}$ is the locally free sheaf
generated by the section with the weight $\mu^{b_0}$.
The cokernel of the complex
$Hom(E_{i}^{\alpha},\nbigg_{\alpha,0})
 \lrarr Hom(E_{i}^{\alpha-1},\nbigg_{\alpha,0})$
is torsion sheaf of length $\eta_{\alpha-1}(i)$.
It is generated by the section with the weight
$\lambda^{-i}\mu^{b_0-b_i-\bar{h}_{\alpha-1,i}}$.
Thus we are done.
\hfill\qed

Thus the contribution of
$\sum_{\alpha=\alpha_-}^0\sum_{i\geq 1}
 \bigl(
Hom(
 E_{i}^{\alpha},\nbigg_{\alpha,0})
 -
 Hom(
 E_{i}^{\alpha-1},
     \nbigg_{\alpha,0})
 \bigr)$
is as follows:
\[
 -\sum_{\alpha= \alpha_-}^0
 \sum_{i\geq 1}\sum_{a=0}^{\eta_{\alpha-1}(i)-1}
 \lambda^{-i}\mu^{b_0-b_i+a-\bar{h}_{\alpha-1,i}}.
\]
When $\alpha>0$,
the sheaf $\nbigg_{\alpha,0}$ is torsion.
Thus the weight decomposition of
$\sum_{i\geq 1}Hom(E_{i}^{\alpha},\nbigg_{\beta,0})$
is similar to the case $i\geq 1,j\geq 1$.
Hence the contribution of
$\sum_{\alpha>0}\sum_{i\geq 1}
 \bigl(
 Hom(
 E_{i}^{\alpha},\nbigg_{\alpha,0})
- Hom(
 E_{i}^{\alpha-1},\nbigg_{\alpha,0})
\bigr)$
is as follows:
\[
\sum_{\alpha>0}
  \sum_{i\geq 1}\sum_{a=0}^{\bar{h}_{\alpha,0}-1}
 \lambda^{-i}
 \bigl(
\mu^{-b_i+b_0-\bar{h}_{\alpha,i}+a}-
\mu^{-b_i+b_0-\bar{h}_{\alpha-1,i}+a}
 \bigr).
\]

We see the case $i=0$.
Note that $E^{\alpha}_0=0$ if $\alpha\leq 0$.
When $\alpha>0$,
the $E_0^{\alpha}$ is locally free sheaf
generated by the section with the weight $\mu^{b_0+\bar{h}_{\alpha,0}}$.
Thus the weight decomposition of
$\sum_{\alpha>0}\sum_{j\geq 0}
Hom(E_{0}^{\alpha},\nbigg_{\alpha,j})$
is as follows:
\[
 \sum_{\alpha>0}
 \sum_{j\geq 0}
 \sum_{a=0}^{\bar{h}_{\alpha,j}-1}
 \lambda^j\mu^{b_j-b_0-\bar{h}_{\alpha,0}+a}.
\]
Since $E_0^{0}=0$,
the weight decomposition of
$\sum_{\alpha>0}\sum_{j\geq 0}
Hom(E_{0}^{\alpha-1},\nbigg_{\alpha,j})$
is as follows:
\[
 \sum_{\alpha>1}
 \sum_{j\geq 0}
 \sum_{a=0}^{\bar{h}_{\alpha,j}-1}
 \lambda^j\mu^{b_j-b_0-\bar{h}_{\alpha-1,0}+a}.
\]
By adding all,
we completed the proof of Proposition \ref{prop;12.3.15}
\hfill\qed

\subsection{The number of cells}

We obtain the following corollary.
\begin{cor}\label{cor;cell_decomposition}

Let $\vecv$ be an element of $\nbiga$.
The punctual parabolic Hilbert schemes
$X^{[\vecv]}_p$
have the cell decomposition.
We have the bijective correspondence
of the set of cells and $\mM(\vecv)$.
The dimension of the cell corresponding to
$\eta\in \mM(\vecv)$ is
$n-||\eta||$.
\end{cor}
\pf
We only calculate the dimension of the cell
corresponding to the fixed point $FP(\eta)$
for each $\eta\in \mM(\vecv)$.
We can use the same discussion as that in \cite{es}.

For the parabolic ideal $(I,\nbigg_{\ast})$
corresponding to the fixed point $FP(\eta)$,
we have the weight decomposition of the tangent space
which is the Laurent polynomial
with variables $\lambda$ and $\mu$,
i.e.,
it is of the form $T_{FP(\eta)}=\sum a_{i\,j}(\eta)\lambda^i\mu^j$.
Note that $a_{i\,j}=0$ for any $(i,j)$
satisfying $i=0$ and $j>0$.
We put 
$S:=\{(i,j)\,|\,\exists \eta\in \mM(\vecv),a_{i\,j}(\eta)\neq 0\}$.
The set $S$ is finite.
For any pair of integers $(w_1,w_2)$,
we have the $G_m$-action on $A^2$
given by
$t\cdot(x,y)=(t^{w_2}x,t^{w_1}y)$.
We can take $w_1,w_2$
such that the following condition holds:
For an element $(i,j)\in S$,
the inequality $i\cdot w_2+j\cdot  w_1>0$ holds
if and only if  $i>0$.

Then we have the induced $G_m$-action on the tangent space
$T_{FP(\eta)}$.
We count the number of the positive weight
by using the isomorphism $T_{FP(\eta)}\simeq \hyperh^0(C(I_{\ast}))$.
The contribution of $\hyperh^0(C_1)$
was calculated by Ellingsrud and Str\o{}mme \cite{es};
the result is $\sum_{i>0}b_i$.
As for the contribution of $\hyperh^{\cdot}(C_{par})$,
it is easy to see that 
the number of the positive weights is the following:
\[
  \sum_{j\geq 1}
 \sum_{a=0}^{\bar{h}_{1,j}-1}
 1=\sum_{j\geq 1}\bar{h}_{1,j}
=\sum_{j\geq 1}\sum_{\alpha>0}\eta_{\alpha}(j).
\]
Thus the number of the positive weights in the tangent space
$T_{FP(\eta)}$ is the following:
\begin{multline}
 \sum_{i>0}b_i+\sum_{j\geq 1}\sum_{\alpha>0}
 \eta_{\alpha}(j)=
 \sum_{j>0}(j-1)\sum_{\alpha} \eta_{\alpha}(j)
 +
 \sum_{j\geq 1}\sum_{\alpha>0} \eta_{\alpha}(j)\\
 =
 \sum_{j>0}j\sum_{\alpha} \eta_{\alpha}(j)
 -\sum_{\alpha\leq 0}\sum_{j\geq 1} \eta_{\alpha}(j)
=
 n-\sum_{\alpha\leq 0}\sum_{j\geq 1} \eta_{\alpha}(j).
\end{multline}
The right hand side can be rewritten as
$n-||\eta||$.
Hence we are done.
\hfill\qed

\vspace{.1in}
For any element $\vecv$ of $\TT$,
the number $\epsilon_{\vecv}$ is defined to be
$0$ $(\vecv\neq 0)$ or $1$ $(\vecv=0)$.
\begin{cor} \label{cor;7.4.1}
Assume that $\rho_-(\vecv)=0$.
\begin{enumerate}
\item
In the decomposition of $X^{[\vecv]}$ given above,
the dimensions of the cells are less than
$n-\epsilon_{\rho_+(\vecv)}+\epsilon_{\vecv}$.
\item
When $\rho_+(\vecv)\neq 0$,
the dimension of the cell corresponding to $\eta\in \mM(\vecv)$
is $n-\epsilon_{\rho_+(\vecv)}+\epsilon_{\vecv}$
if and only if $\eta_{0}=0$.
\item
When $\rho_+(\vecv)=0$ and $\vecv\neq 0$,
the dimension of the cell corresponding to $\eta\in\mM(\vecv)$
is $n-1$ if and only if
$\eta_{0}(i)$ is $0$ ($i\neq n$) or $1$ $(i=n)$.
\hfill\qed
\end{enumerate}
\end{cor}

The cells which have the maximal dimension are called 
the top cells.

\begin{cor}
If $\vecv$ is $n\cdot \vece_0$ or
$n\cdot\vece_0+\vece_{\alpha}\,\,(\alpha>0)$,
then there exists the unique top cell in the decomposition
of the punctual Hilbert scheme $X^{[\vecv]}_p$.
\hfill\qed
\end{cor}


\subsection{Brian\c{c}on's theorem}

Let $\vecv$ be an element of $\nbiga$.
\begin{prop}
The union of the top cells is dense in 
the punctual parabolic Hilbert scheme $X^{[\vecv]}_p$.
\end{prop}
\pf
If $\rho(\vecv)=0$,
the claim is the theorem of Brian\c{c}on.
So we only have to consider the case $\rho(\vecv)\neq 0$.
Due to the isomorphism given in subsection \ref{subsection;7.5.1},
we only have to consider the case $\rho_-(\vecv)=0$
(See Remark \ref{rem;7.5.2}).
We divide the argument into the following steps:
(i) The case $\vecv=n\cdot\vece_0+\vece_{\alpha}$,
(ii) The case $\vecv=n\cdot\vece_0+\sum_{\alpha=1}^{\alpha_+}\vece_{\alpha}$,
(iii) The general case.

\noindent
(i)
If $\vecv=n\vece_0+\vece_{\alpha}$,
then $X^{[\vecv]}_p$ is isomorphic to the punctual nested Hilbert
scheme.
Thus the irreducibility was proved by Ellingsrud-Str\o{}mme
\cite{es2}.

\noindent
(ii)
For $\vecv=m\cdot\vece_0+\sum_{\alpha=1}^{\alpha_+}\vece_{\alpha}$,
we put $\vecv_0:=n\vece_0+\vece_{\alpha_+}$.
We have the universal ideal sheaf $I$
defined over $X_p^{[\vecv_0]}\times X$,
and the universal quotient
$I\otimes\nbigo_D\lrarr \nbigg_{\alpha_+}$
defined over $X_p^{[\vecv_0]}\times D$.
We put
$K:=\ker(I\otimes\nbigo_D\lrarr\nbigg_{\alpha_+})$.
On $X_p^{[\vecv_0]}\times D$,
we take the resolution $(F\lrarr V)$ of $K$,
that is, $V/F\simeq K$.
We consider the moduli scheme $Filt$
of the tuple
$(U,\nbigg_{\ast})$
of a $X^{[\vecv_0]}$-scheme $U$
and the filtration
$V|_{U\times D}\rarr\nbigf_1\rarr\nbigf_{2}\rarr\cdots\rarr
 \nbigf_{\alpha_+-1}\rarr\nbigf_{\alpha_+}=0$
over $U\times D$
satisfying the following:
\begin{itemize}
\item $\nbigf_i$ are flat over $U$.
\item
$\length(\ker(\nbigf_{i}\rarr\nbigf_{i+1})|_{\{u\}\times D})=1$
for any $i=1,\ldots,\alpha_+-1$
and for any closed point $u$ of $U$.
\item
The support of $\ker(\nbigf_i\rarr\nbigf_{i+1})|_{\{u\}\times D}$
is $\{x\}$ for any $i=1,\ldots,\alpha_{+}-1$
and for any closed point $u$ of $U$.
\end{itemize}
Then it is easy to see that $Filt$ is smooth scheme
over $X_p^{[\vecv_0]}$.
The relative dimension of $Filt$ over $X_p^{[\vecv_0]}$
is $(\rank(V)-1)\times (\alpha_+-1)$.
We have the natural morphism
$f:Filt\times D\lrarr X_p^{[\vecv_0]}\times D$.
We denote the projection of $Filt\times D$ to $Filt$
by $p_1$.
On $Filt\times D$,
we have the universal filtration $\nbigf_{\ast}$
of $f^{\ast}(V)$.
Thus we obtain the composite of the morphisms
$f^{\ast}F\lrarr f^{\ast}V\lrarr 
 \nbigf_1$.
Hence we obtain the vector bundle
$p_{1\ast}(\nhom(f^{\ast}(F),\nbigf_1))$
and the section $s$.
The rank of the vector bundle
$p_{1\ast}(\nhom(f^{\ast}(F),\nbigf_1))$
is $\rank(F)\cdot (\alpha_+-1)=(\rank(V)-1)\cdot(\alpha_+-1)$.
Thus the dimensions of
any irreducible components of the $0$-set $s^{-1}(0)$
are larger than $\dim(X_p^{[\vecv_0]})$.

It is easy to see that $s^{-1}(0)$ is isomorphic
to $X_p^{[\vecv]}$.
Thus the dimension of
each irreducible component of $X_p^{[\vecv]}$
is larger than $\dim(X_p^{[\vecv_0]})=n$.
On the other hand,
we have the cell decomposition of $X_p^{[\vecv]}$.
The dimensions of top cells are $n$,
and the dimensions of the other cells are strictly less than $n$.
Thus we can conclude that
the $X_p^{[\vecv]}$ is the closure of the top cells.

\noindent
(iii)
For a general element $\vecv$ of $\nbiga$,
we put $\alpha_+=\sum_{\alpha>0}\rho_{\alpha}(\vecv)$,
and 
$\vecv'=
 \rho_0(\vecv)\cdot \vece_{0}+
 \sum_{\alpha=1}^{\alpha_+}\vece_{\alpha}$.
It is easy to see that we can take a surjective morphism
of $X^{[\vecv']}$ onto $X^{[\vecv]}$
such that the following is satisfied:
\begin{quote}
It is equivariant,
and any fixed point corresponding to a top cell
of $X^{[\vecv']}$
is mapped to a fixed point corresponding to a top cell
of $X^{[\vecv]}$.
\end{quote}
It follows that the points contained in the top cells  are mapped
to the points contained in the top cells.
Since $X^{[\vecv']}$ is the closure of the top cells,
we can conclude that
$X^{[\vecv]}$ is the closure of the top cells
of $X^{[\vecv]}$.
\hfill\qed


\subsection{The class in the Grothendieck ring $K_0(V_k)$}

\label{section;7.1.11}

Let $A^{l}$ denote an $l$-dimensional affine space over $k$.
For a variety $Y$ and an integer $a$,
we put as follows:
\[
 S_a(Y,t)=\sum_{n=0}^{\infty} (Y\times A)^{(n)}\cdot t^n
 \in K_0(V_k)[[t]].
\]
Let $X$ be a projective smooth surface
and $D$ be a divisor.
We put as follows:
\[
 H_0:=
 \sum_{\vecv\in \nbiga}
 X^{[\vecv]} \cdot x^{\vecv}
=\sum_{\vecv\in\nN}X^{[\vecv]}\cdot x^{\vecv}
\in K_0(V_k)[[\nbigp]].
\]
Here we put $X^{[\vecv]}:=\emptyset$ for any $\vecv\not\in\nbiga$.
\begin{thm} \label{thm;12.3.10}
We have the following equality
in  $K_0(V_k)[[\nN]]$:
\[
H_0=
 \prod_{i\geq 1}S_{i-1}(X,x_0^i)
\times
 \prod_{\alpha<0}\prod_{i\geq 1}S_{i-1}(D,x_0^i\cdot x_{\alpha})
\times
 \prod_{\alpha>0}\prod_{i\geq 0}S_{i}(D,x_0^i\cdot x_{\alpha}).
\]
In particular,
we obtain the formula of the $E$-polynomials of
the parabolic Hilbert schemes $X^{[\vecv]}$
to describe
in terms of those of the symmetric powers of $X$ and $D$.
\end{thm}
\pf
Let $\vecv$ be an element of $\nN$.
From the decomposition $X=U+D$,
we have the following direct sum decomposition:
\[
 X^{(\rho_0(\vecv))}\times D^{(\rho(\vecv))}=
 \coprod_{a+b=\rho_0(\vecv)}
 U^{(a)}\times D^{(b)}\times D^{(\rho(\vecv))}=
\coprod_{\substack{a+\rho_0(\vecv')=\rho_0(\vecv),
 \\ \rho(\vecv')=\rho(\vecv)}}
 U^{(a)}\times D^{(\vecv')}
\]
The direct sum decompositions of
$U^{(a)}\times D^{(\vecv')}$
give the following decomposition of $X^{(n)}\times D^{(l_{\ast})}$:
\[
 X^{(\rho_0(\vecv))}\times D^{(\rho(\vecv))}
=
 \coprod_{
\substack{a+\rho_0(\vecv')=\rho_0(\vecv)\\
 (f,g)\in \nbigs(\nnN',a)\times \nbigs(\nN',\vecv') }}
 U^{(f)}_{\ast}\times D_{\ast}^{(g)}.
\]
We have the Hilbert Chow morphism
$\omega:X^{[\vecv]}
\lrarr X^{(\rho_0(\vecv))}\times D^{(\rho(\vecv))}$.
The inverse image
$\omega^{-1}(U^{(f)}_{\ast}\times D^{(g)}_{\ast})$
with the reduced structure
is denoted by $X[f,g]$.
Then we obtain the following equality:
\begin{equation}
H_0=
 \sum_{(f,g)\in \nbigs(\nnN')\times\nbigs(\nN')}
 X[f,g]\cdot
 x_0^{\phi_{\{0\}}(f)}\cdot
 x^{\phi_{\LLA}(g)}.
\end{equation}
Let see $X[f,g]$ in the following two cases:
\begin{itemize}
\item
$a=n$, $b=0$, $f(n)=1$ and $f(a')=0$ for $a'\neq n$.
We denote such $f$ by $\delta_n$.
In this case $g=0$.
\item
$a=0$, $b=n$, $g(\vecv)=1$ and $g(\vecv')=0$ for
$\vecv'\neq \vecv$.
We denote such $g$ by $\delta_{\vecv}$.
In this case $f=0$.
\end{itemize}

We have the fibrations
$X[\delta_n,0]\lrarr U$ and $X[0,\delta_{\vecv}]\lrarr D$.
The fibers are the punctual Hilbert scheme and
the punctual parabolic Hilbert scheme respectively.
The statement for $X[\delta_n,0]$
in the following lemma was proved by  G\"ottsche 
\cite{g}.
\begin{lem} \label{lem;6.30.1}
The fibrations $X[\delta_n,0]\lrarr U$ and
$X[0,\delta_{\vecv}]\lrarr D$ are Zariski trivial.
\end{lem}
\pf
Since the ground field $k$ is algebraically closed,
we can take a Zariski covering $\{U_i\}$
such that we have an etale map $\psi_i:(U_i,U_i\cap D)\lrarr (A^2,A^1)$
for each $i$.
Since the families $X[\delta_n,0]|_{U_i\cap U}\lrarr U\cap U_i$
and $X[0,\delta_{\vecv}]|_{U_i\cap D}\lrarr D\cap U_i$
are the pull backs of the families
$A^2[\delta_n,0]\lrarr A^2$ and $A^2[0,\delta_{\vecv}]\lrarr A^1$
via the morphism $\psi_i$,
we only have to prove that the families
$X[\delta_n,0]\lrarr U$ and $X[0,\delta_{\vecv}]\lrarr D$
are trivial when $X=A^2$ and $D=A^1$.
We can show such trivialities by using the
translations $A^2\lrarr A^2,x\longmapsto x+a,\,(a\in A^2)$
and $A^1\lrarr A^1,x\longmapsto x+a,\,(a\in A^1)$.
\hfill\qed

Let $R(i)$ denote the punctual Hilbert scheme
of the $0$-schemes of length $i$,
and $R(\vecv)$ denote
the punctual parabolic Hilbert scheme of the parabolic length $\vecv$.
Due to Lemma \ref{lem;6.30.1},
we obtain the following equalities in $K_0(V_k)$:
\[
 X[\delta_i,0]= U\times R(i),
\quad
 X[0,\delta_{\vecv}]=D\times R(\vecv).
\]
Due to the result of Ellingsrud-Str\o{}mme \cite{es}
and Proposition \ref{cor;cell_decomposition},
we have the following cell decompositions:
\[
 R(i)=\coprod_{\rho\in \mmM(i)}A^{\phi(\rho)-|\rho|},
\quad
 R(\chi)=
\coprod_{\eta\in \mM(\vecv)}A^{\rho_0(\Psi(\eta))-||\eta||}.
\]
Here we put
$||\eta||:=\sum_{\alpha\leq 0}|\eta_{\alpha}|$.
We have the following isomorphism of the varieties:
\[
 X[f,g]\simeq
 \Bigl(\prod_{i\in \nnN'} X[\delta_i,0]^{(f(i))}\Bigr)_{\ast}
\times
 \Bigl(\prod_{\vecv\in \nN'}X[0,\delta_{\vecv}]^{(g(\vecv))} \Bigr)_{\ast}.
\]
Thus we have the following equality in $K_0(V_k)$:
\begin{equation} \label{eq;12.2.1}
 X[f,g]=
 \Bigl(\prod_{i\in \nnN'} (U\times R_i)^{(f(i))}\Bigr)_{\ast}\times
 \Bigl(\prod_{\vecv\in \nN'}(D\times R(\vecv) )^{(g(\vecv))}\Bigr)_{\ast}.
\end{equation}
We put as follows:
\[
 H_1:=\!\!\!\!\!
 \sum_{f\in \nbigs(\nnN')}
 \Bigl[
 \prod_{i\in \nnN'}
 (U\times R(i))^{(f(i))}
 \Bigr]_{\ast}
 \cdot x_0^{\phi_{\{0\}}(f)},
\quad
 H_2:=\!\!\!\!\!
 \sum_{g\in \nbigs(\nN')}
 \Bigl[
 \prod_{\vecv\in \nN'}
 (D\times R(\vecv))^{(g(\vecv))}
 \Bigr]_{\ast}\cdot x^{\phi_{\LLA}(g)}.
\]
Then we have the equality $H_0=H_1\cdot H_2$.
By applying Lemma \ref{lem;6.30.20},
we shall obtain the formula for $H_1$ and $H_2$.
The generating function $H_1$ is of the following form:
\[
 H_1=\sum_{f\in\nbigs(\nnN')}
  \Bigl[
  \prod_{i\in \nnN'}\Bigl(
 \sum_{\rho\in \mmM(i)}U\times A^{\phi_{\{0\}}(\rho)-|\rho|}
 \Bigr)^{(f(i))}
  \Bigr]_{\ast}
 \cdot 
 x_0^{\phi_{\{0\}}(\rho)}.
\]
Recall that $\mmM=\nbigs(\nbigt(\{0\}))$.
For any element $a\in \nbigt(\{0\})\simeq\seisuu_{\geq\,1}$,
we have the equalities
$\rho(\delta_{a})-|\delta_a|=a-1$,
and $\phi_{\{0\}}(\delta_a)=a$.
Thus we obtain the following formula:
\[
 H_1=
 \prod_{a\in \seisuu_{\geq\,1}}
 \Bigl(
 \sum_{j=0}^{\infty}
 (U\times A^{a-1})^{(j)}t^{aj}
 \Bigr)=\prod_{a\geq 1}S_{a-1}(U,t^a).
\]
The generating function $H_2$ is of the following form:
\[
  H_2=
 \sum_{g\in \nbigs(\nN')}
 \Bigl[
 \prod_{\vecv\in \nN'}
 \Bigl(
 \sum_{\eta\in \mM(\vecv)}D\times A^{\rho_0(\Psi(\eta))-||\eta||}
 \Bigr)^{(g(\vecv))}
 \Bigr]_{\ast}\cdot x^{\phi_{\LLA}(g)}.
\]
By applying Lemma \ref{lem;6.30.20},
we obtain the following formula:
\[
 H_2=
 \prod_{\vecu\in \LLB}
 \Bigl(
 \sum_{j=0}^{\infty}
 \bigl(D\times A^{\rho_0(\Psi(\delta_{\vecu}))-||\delta_{\vecu}||}
 \bigr)^{(j)}\cdot
 x^{\Psi(\delta_{\vecu})}
 \Bigr).
\]
We have the decomposition
$ \LLB= \coprod_{\alpha\in\seisuu}C_{\alpha}$.
For any element $\vecu\in C_{\alpha}$,
we have the equality $\Psi(\delta_{\vecu})=\vecu$.
Hence we have $x^{\Psi(\delta_{\vecu})}=
x_0^{\rho_0(\vecu)}\cdot x_{\alpha}$
and $\rho_0(\Psi(\delta_{\vecu}))=\rho_0(\vecu)$
for $\vecu\in C_{\alpha}$.
We also have the following equality for $\vecu\in C_{\alpha}$:
\[
 ||\delta_{\vecu}||=
\left\{
 \begin{array}{ll}
 0 &(\alpha>0) \\
 1 &(\alpha\leq 0).
 \end{array}
\right.
 \]
Thus we obtain the following formula:
\begin{equation}
 H_2=
 \prod_{i\geq 1}  S_{i-1}(D,x_0^i)
\times
 \prod_{\alpha<0}
 \prod_{i\geq 1} S_{i-1}(D,x_0^i\cdot x_{\alpha})
\times
 \prod_{\alpha>0}
 \prod_{i\geq 0}S_{i}(D,x_0^i\cdot x_{\alpha}).
\end{equation}
Note that we have the following equality:
\[
 S_{a}(U,t)\cdot S_{a}(D,t)=
 S_a(U+D,t)=S_a(X,t).
\]
In all we obtain the formula desired.
\hfill\qed


\section{The incidence varieties and their composition}

In this section,
we only consider the parabolic Hilbert schemes $X^{[\vecv]}$
such that $\rho_-(\vecv)=0$.
We put $\TTt:=\TT_0\oplus \TT_+$,
which is a subgroup of $\TT$.
We put $C_L:=C_0\cup C_+=C\cap \TTt$,
$-C_L:=(-C)\cap \TTt$
and $\widetilde{C}_L=C_L\cup (-C_L)$.
We put $\nbiga_L:=\nbiga\cap \TTt$.

For any element $\vecv$ of $\TTt$,
the number $\epsilon_{\vecv}$ is defined to be
$1$ $(\vecv=0)$ or $0$ $(\vecv\neq 0)$.

\subsection{The incidence varieties}

\subsubsection{$\vecu\in C_L$}
Let $\vecv$ be an element of $\nbiga_L$.
We put
$\nbiga_L(\vecv):=\{\veca\in \nbiga_L\,|\,\vecv-\veca\in \nbiga_L\}$.
For any $\vecu\in C_L$
and any $\veca\in \nbiga_L(\vecv)$, 
the locally closed subscheme
$\widetilde{Z}(\vecv,\vecu)_{\veca}$
of
$X^{[\vecv]}\times Y_{\vecu}\times X^{[\vecv+\vecu]}$
is defined as follows:
\[
\widetilde{Z}(\vecv,\vecu)_{\veca}
:=
 \left\{
  \bigl(
   (I_1,\nbigg_{1\,\ast}),x,(I_2,\nbigg_{2\,\ast})
  \bigr)
 \,\Bigl|\Bigr.\,
 \begin{array}{l}
   (I_1,\nbigg_{1,\ast})=(I_2,\nbigg_{2,\ast})
 \mbox{ on $X-\{x\}$, } \\
 \length_x(\nbigo/I_1,\nbigg_{1,\ast})=\veca,\,
 \length_x(\nbigo/I_2,\nbigg_{2,\ast})=\veca+\vecu
 \end{array}
 \right\}.
\]

\begin{lem} \label{lem;3.1.1}
Let $\vecv$, $\vecu$ and $\veca$ be as above.
The dimension of $\widetilde{Z}(\vecv,\vecu)_{\veca}$ is as follows:
\[
d(\vecv)+\rho_0(\vecu)-|\veca|_+ +1
 -\epsilon_{\rho_+(\veca)}+\epsilon_{\veca}
 -\epsilon_{\rho_+(\veca+\vecu)}
 +\epsilon_{\rho_+(\veca)}\cdot\epsilon_{\rho_+(\veca+\vecu)}.
\]
Let $\vecv$ be an element of $\nbiga_L$
and $\vecu$ be an element of $C_L$.
It holds that
$\dim(\widetilde{Z}(\vecv,\vecu)_0)=d(\vecv)+\rho_0(\vecu)+1$.
Let $\veca\neq 0$ be an element of $\nbiga_L(\vecv)$.
Then we have the inequality
$\dim(\widetilde{Z}(\vecv,\vecu)_{\veca})<d(\vecv)+\rho_0(\vecu)+1$.
\end{lem}
\pf
It is easy to see 
the following equality:
\[
 \dim(\widetilde{Z}(\vecv,\vecu)_{\veca})=
 \dim(X^{[\vecv-\veca]})
+\dim(X_p^{[\veca]})+\dim(X_p^{[\veca+\vecu]})
+1+\epsilon_{\rho_+(\veca)}\cdot\epsilon_{\rho_+(\vecu+\veca)}.
\]
Due to Corollary \ref{cor;7.4.1},
the right hand side is equal to the following:
\begin{multline}
 2\rho_0(\vecv-\veca)
+\sum_{\alpha> 0}\rho_{\alpha}(\vecv-\veca)
+\rho_0(\veca)-\epsilon_{\rho_+(\veca)}+\epsilon_{\veca}
+\rho_0(\vecu+\veca)-\epsilon_{\rho_+(\vecu+\veca)}
+\epsilon_{\vecu+\veca}
+1+\epsilon_{\rho_+(\veca+\vecu)}\cdot\epsilon_{\rho_+(\veca)}\\
= d(\vecv)+\rho_0(\vecu)+1
-|\veca|_+
 -\epsilon_{\rho_+(\veca)}+\epsilon_{\veca}
 -\epsilon_{\rho_+(\vecu+\veca)}+\epsilon_{\vecu+\veca}
 +\epsilon_{\rho_+(\veca+\vecu)}\cdot\epsilon_{\rho_+(\veca)}.
\end{multline}
Since the elements $\vecu$ and $\veca$ are contained in $\nbiga_L$,
the sum $\veca+\vecu$ is not $0$.
Thus we have $\epsilon_{\veca+\vecu}=0$.
Hence we obtain the first equality.
It implies the equality in the case $\veca=0$.

If $\veca\neq 0$, then $\epsilon_{\veca}=0$.
We have the inequalities $-|\veca|_+-\epsilon_{\rho_+(\veca)}\leq -1$
and $-\epsilon_{\rho_+(\veca+\vecu)}\cdot(1-\epsilon_{\rho_+(\veca)})\leq 0$.
Thus we obtain the inequality in the case $\veca\neq 0$.
\hfill\qed

\vspace{.1in}
We denote the closure of $\widetilde{Z}(\vecv,\vecu)_0$
by $Z(\vecv,\vecu)$,
whose dimension is $d(\vecv)+\rho_0(\vecu)+1$.
We call it the incidence variety
for $(\vecv,\vecu)$
in the case $\vecu\in C$.

\subsubsection{$\vecu\in -C_L$}
Let $\vecv$ be an element of $\nbiga_L$
and $\vecu$ be an element of $-C_L$
such that $\vecv+\vecu$ is contained in $\nbiga_L$.
The incidence variety
$Z(\vecv,\vecu)\subset X^{[\vecv]}\times Y_{\vecu}\times
  X^{[\vecv+\vecu]}$
is defined as follows:
We have the incidence variety
$Z(\vecv+\vecu,-\vecu)$ which is a closed subvariety
of $X^{[\vecv+\vecu]}\times Y_{-\vecu}\times X^{[\vecv]}$.
We have the isomorphism
of $X^{[\vecv+\vecu]}\times Y_{\vecu}\times X^{[\vecv]}$
to $X^{[\vecv]}\times Y_{\vecu}\times X^{[\vecv+\vecu]}$
by the transposition of the first and third components.
The image of $Z(\vecv+\vecu,-\vecu)$
is denoted by $Z(\vecv,\vecu)$,
which we call the incidence variety
in the case $\vecu\in -C_L$.

\subsubsection{The dimension of the incidence variety}

Let $\vecv$ be an element of $\nbiga_L$.
We put
$\widetilde{C}_L(\vecv):=\{\vecu\in \widetilde{C}_L\,|
  \,\vecv+\vecu\in \nbiga_L\}$.
Recall that 
the function $g:\widetilde{C}_L\lrarr \{-1,0\}$
is defined as follows:
\[
 g(\vecu):=
\left\{
 \begin{array}{rl}
  -1 & (\vecu\in -C_L,\,\rho(\vecu)\neq 0),\\
  0 & (\mbox{otherwise}).
 \end{array}
\right.
\]
\begin{lem} \label{lem;7.8.15}
Let $\vecv$ be an element of $\nbiga$,
and $\vecu$ be an element of $\widetilde{C}(\vecv)$.
The dimension of the incidence variety $Z(\vecv,\vecu)$
is $d(\vecv)+\rho_0(\vecu)+1+g(\vecu)$.
\end{lem}
\pf
In the case $\vecu\in C_L$,
the equality is proved in Lemma \ref{lem;3.1.1}.
Let us consider the case $\vecu=-\vecu'\in -C_L$.
The dimension is same as the following:
\begin{multline}
 \dim(Z(\vecv-\vecu',\vecu'))=
 d(\vecv-\vecu')+\rho_0(\vecu')+1
=d(\vecv)-2\rho_0(\vecu')-|\vecu'|_+  +\rho_0(\vecu')+1
\\=
 d(\vecv)+2\rho_0(\vecu)-|\vecu'|_+    
 -\rho_0(\vecu)+1
=d(\vecv)+\rho_0(\vecu)+1-|\vecu'|_+.
\end{multline}
It is easy to see that $g(\vecu)=-|\vecu'|_+$
when $\vecu\in -C_L$.
\hfill\qed

\vspace{.1in}

\noindent
{\bf Notation:}
We put $d(\vecv,\vecu):=d(\vecv)+\rho_0(\vecu)+1+g(\vecu)$
for any $\veca\in \nbiga_L$ and $\vecu\in \widetilde{C}_L(\vecv)$.
\hfill\qed

\subsection{Composition varieties}
\label{subsection;7.8.20}

Let $\vecv$ be an element of $\nbiga_L$.
We have the incidence varieties
$Z(\vecv,\vecu)$ for any
 $\vecv\in \nbiga_L$ and $\vecu\in \widetilde{C}_L(\vecv)$.
We are interested in the correspondence morphisms
induced by the incidence varieties
and the composition of them.
For any $\vecv\in \nbiga_L$,
we put $
 \widetilde{C}^2_L(\vecv):=
 \bigl\{(\vecu_1,\vecu_2)\in\widetilde{C}_L^2\,|\,
 \vecv+\vecu_1\in \nbiga_L,\vecv+\vecu_1+\vecu_2\in \nbiga_L
 \bigr\}$.
Moreover we put as follows:
\[
 Comp_L:=\{(\vecv,\vecu_1,\vecu_2)\,|\,
 \vecv\in \nbiga_L,(\vecu_1,\vecu_2)\in\widetilde{C}^2_L(\vecv) \}.
\]
For any $(\vecv,\vecu_1,\vecu_2)\in Comp_L$,
we consider the product
$X^{[\vecv]}\times Y_{\vecu_1}\times
 X^{[\vecv+\vecu_1]}\times
 Y_{\vecu_2}\times X^{[\vecv+\vecu_1+\vecu_2]}$.
We denote the projections
to $X^{[\vecv]}\times Y_{\vecu_1}\times X^{[\vecv+\vecu_1]}$
and
$X^{[\vecv+\vecu_1]}\times
 Y_{\vecu_2}\times X^{[\vecv+\vecu_1+\vecu_2]}$
by $\pi_{4,5}$ and $\pi_{1,2}$ respectively.
Then we have the following subvariety:
\[
 Z(\vecv,\vecu_1,\vecu_2):=
 \pi_{4,5}^{-1}\bigl(Z(\vecv,\vecu_1)\bigr)
\cap
 \pi_{1,2}^{-1}\bigl(Z(\vecv+\vecu_1,\vecu_2)\bigr).
\]
The expected dimension $d(\vecv,\vecu_1,\vecu_2)$
of $Z(\vecv,\vecu_1,\vecu_2)$
is  given as follows:
\begin{multline}
d(\vecv,\vecu_1,\vecu_2):=\\
\dim \Bigl(\pi_{4,5}^{-1}\bigl(Z(\vecv,\vecu_1)\bigr)
     \Bigr)+
\dim \Bigl(
\pi_{1,2}^{-1}\bigl(Z(\vecv+\vecu_1,\vecu_2)\bigr)
     \Bigr)
-\dim \bigl(
 X^{[\vecv]}\times Y_{\vecu_1}\times
 X^{[\vecv+\vecu_1]}\times
 Y_{\vecu_2}\times X^{[\vecv+\vecu_1+\vecu_2]}\bigr)\\
=
 \dim(Z(\vecv,\vecu_1))+\dim(Z(\vecv+\vecu_1,\vecu_2))-
 X^{[\vecv+\vecu_1]}\\
=d(\vecv)+\rho_0(\vecu_1)+1+g(\vecu_1)
+d(\vecv+\vecu_1)+\rho_0(\vecu_2)+1+g(\vecu_2)
 -d(\vecv+\vecu_1)\\
=d(\vecv)+\rho_0(\vecu_1+\vecu_2)
 +2+g(\vecu_1)+g(\vecu_2).
\end{multline}
We have the intersection product
(\cite{f}):
\[
 [Z(\vecv,\vecu_1,\vecu_2)]:=
 \pi_{4,5}^{\ast}\bigl([Z(\vecv,\vecu_1)]\bigr)
\cdot
 \pi_{1,2}^{\ast}\bigl([Z(\vecv+\vecu_1,\vecu_2)]\bigr)
\in 
CH_{d(\vecv,\vecu_1,\vecu_2)}
 \bigl(Z(\vecv,\vecu_1,\vecu_2)\bigr).
\]
We denote the projection of
$X^{[\vecv]}\times Y_{\vecu_1}\times X^{[\vecv+\vecu_1]}
 \times Y_{\vecu_2}\times X^{[\vecv+\vecu_1+\vecu_2]}$
to $X^{[\vecv]}\times Y_{\vecu_1}
 \times Y_{\vecu_2}\times X^{[\vecv+\vecu_1+\vecu_2]}$
by $\pi_3$.
Then we obtain 
the class
$\pi_{3\,\ast}\bigl([Z(\vecv,\vecu_1,\vecu_2)]\bigr)\in
CH_{d(\vecv,\vecu_1,\vecu_2)}(\pi_3(Z(\vecv,\vecu_1,\vecu_2)))$.


For any $(\vecv,\vecu_1,\vecu_2)\in Comp_L$,
we put as follows:
\[
\begin{array}{l}
 St_0(\vecv,\vecu_1,\vecu_2):=
 \{(\veca_1,\veca_2)\in \nbiga_L^2\,|\,
 \veca_i+\vecu_i\in \nbiga_L,\,
 \vecv-\veca_1-\veca_2\in \nbiga_L \},\\
 St_1(\vecv,\vecu_1,\vecu_2):=
 \{\veca\in \nbiga_L\,|\,
 \vecv-\veca\in \nbiga_L,\,
 \veca+\vecu_1\in \nbiga_L,\,
 \veca+\vecu_1+\vecu_2\in \nbiga_L
 \}.
\end{array}
\]
We consider the following locally closed subschemes
$\widetilde{Z}(\vecv,\vecu_1,\vecu_2)_{0,(\veca_1,\veca_2)}$
and $\widetilde{Z}(\vecv,\vecu_1,\vecu_2)_{1,\veca}$
of 
$X^{[\vecv]}\times Y_{\vecu_1}\times Y_{\vecu_2}\times 
 X^{[\vecv+\vecu_1+\vecu_2]}$:
For any $(\veca_1,\veca_2)\in St_0(\vecv,\vecu_1,\vecu_2)$,
the subscheme $\widetilde{Z}(\vecv,\vecu_1,\vecu_2)_{0,(\veca_1,\veca_2)}$
is defined as follows:
 \begin{multline}
   \widetilde{Z}(\vecv,\vecu_1,\vecu_2)_{0,(\veca_1,\veca_2)}
:=\\
\left\{
 \bigl((I_1,\nbigg_{1,\ast}),x_1,x_2,
       (I_2,\nbigg_{2,\ast})
 \bigr)\,\Bigl|\Bigr.\,
\begin{array}{l}
 x_1\neq x_2,\,\,
 (I_1,\nbigg_{1,\ast})=(I_2,\nbigg_{2,\ast})
 \mbox{ on $X-\{x_1,x_2\}$, }\\
 \length_{x_i}(\nbigo/I_1,\nbigg_{1\,\ast})=\veca_i,\,
 \length_{x_i}(\nbigo/I_{2},\nbigg_{2\,\ast})
 =\veca_{i}+\vecu_i
 \end{array}
 \right\}.
\end{multline}
For any $\veca\in St_1(\vecv,\vecu_1,\vecu_2)$,
the locally closed subscheme
$\widetilde{Z}(\vecv,\vecu_1,\vecu_2)_{1,\veca}$
is defined as follows:
\begin{multline}
 \widetilde{Z}(\vecv,\vecu_1,\vecu_2)_{1,\veca}
 :=\\
 \left\{
 \bigl((I_1,\nbigg_{1,\ast}),x,x,
       (I_2,\nbigg_{2,\ast})
 \bigr)\,\Bigl|\Bigr.\,
 \begin{array}{l}
 (I_1,\nbigg_{1,\ast})=(I_2,\nbigg_{2,\ast})
 \mbox{ on $X-\{x\}$, }\\
 \length_{x}(\nbigo/I_1,\nbigg_{1\,\ast})=\veca,\,
 \length_{x}(\nbigo/I_2,\nbigg_{2\,\ast})=\vecu_1+\vecu_2+\veca
 \end{array}
 \right\}.
\end{multline}
We denote the closures of
$\widetilde{Z}(\vecv,\vecu_1,\vecu_2)_{0,(\veca_1,\veca_2)}$
and $\widetilde{Z}(\vecv,\vecu_1,\vecu_2)_{1,\veca}$
by $\overline{Z}(\vecv,\vecu_1,\vecu_2)_{0,(\veca_1,\veca_2)}$
and $\overline{Z}(\vecv,\vecu_1,\vecu_2)_{1,\veca}$
respectively.
We put as follows:
\[
 \overline{Z}(\vecv,\vecu_1,\vecu_2):=
 \bigcup_{(\veca_1,\veca_2)\in St_0(\vecv,\vecu_1,\vecu_2)}
 \overline{Z}(\vecv,\vecu_1,\vecu_2)_{0,(\veca_1,\veca_2)}
 \cup
 \bigcup_{\veca\in St_1(\vecv,\vecu_1,\vecu_2)}
  \overline{Z}(\vecv,\vecu_1,\vecu_2)_{1,\veca}.
\]
It is easy to see that
the image
$\pi_{3}(Z(\vecv,\vecu_1,\vecu_2))$
is contained in
$\overline{Z}(\vecv,\vecu_1,\vecu_2)$.

The varieties
$\widetilde{Z}(\vecv,\vecu_1,\vecu_2)_{0,(\veca_1,\veca_2)}$,
$\widetilde{Z}(\vecv,\vecu_1,\vecu_2)_{1,\veca}$,
$\overline{Z}(\vecv,\vecu_1,\vecu_2)_{0,(\veca_1,\veca_2)}$,
$\overline{Z}(\vecv,\vecu_1,\vecu_2)_{1,\veca}$
and $\overline{Z}(\vecv,\vecu_1,\vecu_2) $
are called the composition varieties.
In the next subsection,
we show the following:
\begin{prop}[Corollary\ref{cor;3.10.1}, Corollary \ref{cor;3.10.2}  ]
\label{prop;3.10.1}
Let $(\vecv,\vecu_1,\vecu_2)$ be an element of
$Comp_L$.
\begin{enumerate}
\item
The dimension of $\overline{Z}(\vecv,\vecu_1,\vecu_2)$
is $d(\vecv,\vecu_1,\vecu_2)$.
\item
For $\vecu\in \widetilde{C}_L$,
the element $\kappa(\vecu)\in A_L$
is defined to be $\vecu$ $(\vecu\in -C_L)$ or $0$ $(\vecu\in C_L)$.
Then the component
$\overline{Z}(\vecv,\vecu_1,\vecu_2)
    _{0,(\kappa(\vecu_1),\kappa(\vecu_2))}$
has the expected dimension $d(\vecv,\vecu_1,\vecu_2)$.

\item
When $\vecu_1\in C_L$ and $\vecu_1+\vecu_2=0$,
then the component
$\overline{Z}(\vecv,\vecu_1,-\vecu_1)_{1,(0,0)}$
has the expected dimension $d(\vecv,\vecu_1,\vecu_2)$.
It coincides with the diagonal component
of $X^{[\vecv]}\times Y_{\vecu_1}\times Y_{\vecu_1}\times X^{[\vecv]}
\simeq (X^{[\vecv]}\times Y)^2$.
\item
The dimensions of the other components are strictly
less than the expected dimension $d(\vecv,\vecu_1,\vecu_2)$.
\end{enumerate}
\end{prop}

\subsection{The estimate of the dimensions of the composition varieties}

\subsubsection{$\widetilde{Z}(\vecv,\vecu_1,\vecu_2)_{0,(\veca_1,\veca_2)}$}

Let $(\vecv,\vecu_1,\vecu_2)$ be an element of $Comp_L$.
\begin{lem}
For elements $\veca_1$ and $\veca_2$ of $St_0(\vecv,\vecu_1,\vecu_2)$,
the dimension of 
$\widetilde{Z}(\vecv,\vecu_1,\vecu_2)_{0,(\veca_1,\veca_2)}$ is
as follows:
\[
 d(\vecv)+\rho_0(\vecu_1+\vecu_2)+2
+\sum_{i=1,2}\Bigl(
 -|\veca_i|_+ -\epsilon_{\rho(\veca_i)} +\epsilon_{\veca_i}
 -\epsilon_{\rho(\veca_i+\vecu_i)}+\epsilon_{\veca_i+\vecu_i}
 +\epsilon_{\rho(\vecu_i)}\cdot\epsilon_{\rho(\veca_i)}
   \cdot\epsilon_{\rho(\veca_i+\vecu_i)}
 \Bigr).
\]
\end{lem}
\pf
Let $\bigl((I_1,\nbigg_{1,\ast}),(x_1,x_2),(I_2,\nbigg_{2,\ast}) \bigr)$
be a point contained in 
$\widetilde{Z}(\vecv,\vecu_1,\vecu_2)_{0,(\veca_1,\veca_2)}$.
The contribution of the parts
$(I_1,\nbigg_{1\,\ast})|_{X-\{x_1,x_2\}}
=(I_2,\nbigg_{2\,\ast})|_{X-\{x_1,x_2\}}$
to the dimension of
$Z(\vecv,\vecu_1,\vecu_2)_{0,(\veca_1,\veca_2)}$
is same as
$\dim(X^{[\vecv-\veca_1-\veca_2]})=d(\vecv-\veca_1-\veca_2)$.
Note that $x_i$ is contained in $D$
unless $\rho(\vecu_i)=\rho(\veca_i)=\rho(\vecu_i+\veca_i)=0$.
Thus the contribution of 
$(I_1,\nbigg_{1\,\ast})_{x_i}$, $(I_2,\nbigg_{2\,\ast})_{x_i}$
and $x_i$ to the dimension 
of $\widetilde{Z}(\vecv,\vecu_1,\vecu_2)_{0,(\veca_1,\veca_2)}$
is same as the following:
\begin{multline}
 \dim (X^{[\veca_i]}_{x_i})+\dim(X^{[\veca_i+\vecu_i]}_{x_i})+
 (1+\epsilon_{\rho_+(\vecu_i)}\cdot \epsilon_{\rho(\veca_i)}
    \cdot\epsilon_{\rho_+(\vecu_i+\veca_i)})\\
=\rho_0(\veca_{i})-\epsilon_{\rho_+(\veca_i)}+\epsilon_{\veca_i}
+\rho_0(\veca_{i}+\vecu_{i})
-\epsilon_{\rho_+(\veca_i+\vecu_i)}+\epsilon_{\veca_i+\vecu_i}
+1+\epsilon_{\rho_+(\vecu_i)}\cdot \epsilon_{\rho(\veca_i)}
    \cdot\epsilon_{\rho_+(\vecu_i+\veca_i)}\\
=\rho_0(\vecu_{i})+1+2\rho_0(\veca_{i})
 -\epsilon_{\rho_+(\veca_i)}+\epsilon_{\veca_i}
 -\epsilon_{\rho_+(\veca_i+\vecu_i)}+\epsilon_{\veca_i+\vecu_i}
+\epsilon_{\rho_+(\vecu_i)}\cdot \epsilon_{\rho(\veca_i)}
    \cdot\epsilon_{\rho_+(\vecu_i+\veca_i)}.
\end{multline}
Thus the dimension considered is the following:
\begin{multline}
 d(\vecv-\veca_1-\veca_2)+
\sum_{i=1,2}
 \Bigl(
 \rho_0(\vecu_{i})+1+2\rho_0(\veca_{i})
 -\epsilon_{\rho_+(\veca_i)}+\epsilon_{\veca_i}
 -\epsilon_{\rho_+(\veca_i+\vecu_i)}+\epsilon_{\veca_i+\vecu_i}
+\epsilon_{\rho_+(\vecu_i)}\cdot \epsilon_{\rho(\veca_i)}
    \cdot\epsilon_{\rho_+(\vecu_i+\veca_i)}
 \Bigr)\\
=d(\vecv)+\rho_0(\vecu_{1}+\vecu_{2})+2
+\sum_{i=1,2}
 \Bigl(
  -|\veca|_+
 -\epsilon_{\rho_+(\veca_i)}+\epsilon_{\veca_i}
 -\epsilon_{\rho_+(\veca_i+\vecu_i)}+\epsilon_{\veca_i+\vecu_i}
+\epsilon_{\rho_+(\vecu_i)}\cdot \epsilon_{\rho(\veca_i)}
    \cdot\epsilon_{\rho_+(\vecu_i+\veca_i)}
 \Bigr).
\end{multline}
Hence we are done.
\hfill\qed

We need the following estimate.
\begin{lem}
Let $\vecv,\vecu_i,\veca_i$ be in the previous lemma.
We have the following inequality:
\begin{equation} \label{eq;7.5.10}
 -|\veca_i|_+
 -\epsilon_{\rho_+(\veca_i)}+\epsilon_{\veca_i}
 -\epsilon_{\rho_+(\veca_i+\vecu_i)}+\epsilon_{\veca_i+\vecu_i}
+\epsilon_{\rho_+(\vecu_i)}\cdot \epsilon_{\rho(\veca_i)}
    \cdot\epsilon_{\rho_+(\vecu_i+\veca_i)}
 \leq g(\vecu_i).
\end{equation}
It is an equality if and only if $\veca_i=0$ or $\veca_i+\vecu_i=0$.
\end{lem}
\pf
We divide the claim into the three cases:
(1) $\vecu_i\in C_L$,\,
(2) $\vecu_i\in -C_L$ and $\rho_+(\vecu_i)=0$,\,
(3) $\vecu_i\in -C_L$ and $\rho_+(\vecu_i)\neq 0$.

\noindent
{\bf Case (1)}
We have the equality $g(\vecu_i)=0$.
Since $\veca_i\in \nbiga_L$, the element $\veca_i+\vecu_i$ is not $0$.
Thus $\epsilon_{\veca_i+\vecu_i}=0$.
If $\veca_i\neq 0$,
the left hand side of (\ref{eq;7.5.10}) is as follows:
\[
 -|\veca|_+ -\epsilon_{\rho_+(\veca_i)}-\epsilon_{\rho_+(\veca_i+\vecu_i)}
 +\epsilon_{\rho_+(\vecu_i)}\cdot \epsilon_{\rho(\veca_i)}
    \cdot\epsilon_{\rho_+(\vecu_i+\veca_i)}.
\]
We have the inequalities
$-|\veca|_+ -\epsilon_{\rho_+(\veca_i)}\leq -1$
and $ -\epsilon_{\rho_+(\veca_i+\vecu_i)}
 +\epsilon_{\rho_+(\vecu_i)}\cdot \epsilon_{\rho(\veca_i)}
    \cdot\epsilon_{\rho_+(\vecu_i+\veca_i)}\leq 0$.
Thus the left hand side of (\ref{eq;7.5.10})
is strictly less than $0$.

If $\veca_i=0$, then the left hand side 
of (\ref{eq;7.5.10}) is equal to the following:
$-0-1+1-\epsilon_{\rho_+(\vecu_i)}+\epsilon_{\rho_+(\vecu_i)}^2=0$.

\noindent
{\bf Case (2)}
We have the equality $g(\vecu_i)=0$.
Since $\veca_i+\vecu_i$ is contained in $\nbiga_L$,
the element $\veca_i$ is not $0$. Thus we have $\epsilon_{\veca_i}=0$.

If $\veca_i+\vecu_i\neq 0$,
then the left hand side of (\ref{eq;7.5.10}) is same as the following:
\[
 -|\veca_i|_+ -\epsilon_{\rho_+(\veca_i)}
 -\epsilon_{\rho_+(\veca_i+\vecu_i)} 
 +\epsilon_{\rho_+(\veca_i)}\cdot\epsilon_{\rho_+(\veca_i+\vecu_i)}.
\]
Here we have the inequalities
$-|\veca_i|_+ -\epsilon_{\rho_+(\veca_i)}\leq -1 $
and $-\epsilon_{\rho_+(\veca_i+\vecu_i)}
 +\epsilon_{\rho_+(\veca_i)}\cdot\epsilon_{\rho_+(\veca_i+\vecu_i)}\leq 0$.
Thus the left hand side of (\ref{eq;7.5.10})
is strictly less than $0$.

If $\veca_i+\vecu_i=0$,
then $\rho_+(\veca_i)=0$.
It implies $|\veca_i|_+=0$
and $\epsilon_{\rho_+(\veca_i)}=1$.
Thus the left hand side is same as the following:
$ -0-1+0-1+1+1=0$.

\noindent
{\bf Case (3)}
We have $g(\vecu_i)=-1$.
Since $\vecu_i+\veca_i$ is contained in $\nbiga_L$,
we have $\veca_i\neq 0$ and $\rho_+(\veca_i)\neq 0$.

If $\veca_i+\vecu_i\neq 0$,
the left hand side of (\ref{eq;7.5.10})
is same as 
$ -|\veca_i|_+  -\epsilon_{\rho_+(\veca_i+\vecu_i)}$.
If $|\veca_i|_+\geq 2$,
then the left hand side is strictly less than $-1$.
Assume that $|\veca_i|_+=1$.
Since $\veca_i+\vecu_i\in \nbiga_L$ 
and since we have $|\vecu_i|_+\leq 1$,
$\rho_+(\veca_i+\vecu_i)$ must be $0$
and thus $\epsilon_{\rho_+(\veca_i+\vecu_i)}=1$.
Thus the left hand side is strictly less than $-1$.

If $\veca_i+\vecu_i=0$,
we know the equalities  $|\veca_i|_+=1$
and
$-|\veca_i|_+ -\epsilon_{\rho_+(\veca_i+\vecu_i)}+\epsilon_{\veca_i+\vecu_i}
=-1 -1+1=-1$.
\hfill\qed

\vspace{.1in}
We have the following direct corollary.
\begin{cor}\label{cor;3.10.1}
Let $(\vecv,\vecu_1,\vecu_2)$
be an element of $Comp_L$.
Let $(\veca_1,\veca_2)$ be an element of
$St_0(\vecv,\vecu_1,\vecu_2)$.
We have the following inequality:
\[
 \dim(\widetilde{Z}(\vecv,\vecu_1,\vecu_2)_{0,(\veca_1,\veca_2)} )
\leq d(\vecv,\vecu_1,\vecu_2).
\]
It is an equality if and only if
$\veca_i=0$ or $\veca_i+\vecu_i=0$ for any $i\in\{1,2\}$,
i.e.,
$(\veca_1,\veca_2)=(\kappa(\vecu_1),\kappa(\vecu_2))$.
\hfill\qed
\end{cor}

\vspace{.1in}
\noindent
{\bf Notation:}
The locally closed subscheme
$\widetilde{Z}(\vecv,\vecu_1,\vecu_2)_{0,(\kappa(\vecu_1),\kappa(\vecu_2))}$
is denoted by $\widetilde{Z}(\vecv,\vecu_1,\vecu_2)_0$.
The closure
$\overline{Z}(\vecv,\vecu_1,\vecu_2)_{0,(\kappa(\vecu_1),\kappa(\vecu_2))}$
is denoted by $\overline{Z}(\vecv,\vecu_1,\vecu_2)_0$.
\hfill\qed

\subsubsection{$\widetilde{Z}(\vecv,\vecu_1,\vecu_2)_{1,\veca}$ }

Let $(\vecv,\vecu_1,\vecu_2)$
be an element of $Comp_L$ as usual.

\begin{lem}
Let $\veca$ be an element of $St_1(\vecv,\vecu_1,\vecu_2)$.
Then the dimension of $\widetilde{Z}(\vecv,\vecu_1,\vecu_2)_{1,\veca}$
is as follows:
\[
d(\vecv)+\rho_0(\vecu_1+\vecu_2)+2+
\Bigl(
 -|\veca|_+-1-\epsilon_{\rho_+(\veca)}+\epsilon_{\veca}
 -\epsilon_{\rho_+(\veca+\vecu_1+\vecu_2)}+\epsilon_{\veca+\vecu_1+\vecu_2}
 +\epsilon_{\rho_+(\vecu_1)}\cdot\epsilon_{\rho_+(\vecu_2)}
   \cdot\epsilon_{\rho_+(\veca)}
\Bigr).
\]
\end{lem}
\pf
Let $\bigl((I_1,\nbigg_{1,\ast}),(x,x),(I_2,\nbigg_{2,\ast})\bigr)$
be a point of
$\widetilde{Z}(\vecv,\vecu_1,\vecu_2)_{1,\veca}$.
The contribution
of 
$(I_1,\nbigg_{1,\ast})|_{X-\{x\}}=(I_2,\nbigg_{2,\ast})|_{X-\{x\}}$
is same as $\dim(X^{[\vecv-\veca]})=d(\vecv-\veca)$.
The point $x$ is contained in $D$ unless
$\rho_+(\veca)=\rho_+(\vecu_1)=\rho_+(\vecu_2)
=\rho_+(\veca+\vecu_1+\vecu_2)=0$.
Thus the contribution
of $(I_1,\nbigg_{1,\ast})_x$,
$ (I_2,\nbigg_{2,\ast})_x$ and $x$ to the dimension
of $\widetilde{Z}(\vecv,\vecu_1,\vecu_2)_{1,\veca}$
is same as the following:
\begin{multline}
 \bigl( \rho_0(\veca)-\epsilon_{\rho_+(\veca)}+\epsilon_{\veca}\bigr)
+\bigl(
\rho_0(\veca+\vecu_1+\vecu_2)
-\epsilon_{\rho_+(\veca+\vecu_1+\vecu_2)}
+\epsilon_{\veca+\vecu_1+\vecu_2}\bigr)
+\bigl(
1+\epsilon_{\rho_+(\vecu_1)}\cdot\epsilon_{\rho_+(\vecu_2)}\cdot
   \epsilon_{\rho_+(\veca)}\cdot\epsilon_{\rho_+(\veca+\vecu_1+\vecu_2)}
\bigr)\\
=2\rho_0(\veca)+\rho_0(\vecu_1+\vecu_2)+1
 -\epsilon_{\rho_+(\veca)}+\epsilon_{\veca}
 -\epsilon_{\rho_+(\veca+\vecu_1+\vecu_2)}
 +\epsilon_{\veca+\vecu_1+\vecu_2}
 +\epsilon_{\rho_+(\vecu_1)}\cdot\epsilon_{\rho_+(\vecu_2)}\cdot
   \epsilon_{\rho_+(\veca)}.
\end{multline}
Thus the dimension considered is the following:
\begin{multline}
 d(\vecv-\veca)+
 2\rho_0(\veca)+\rho_0(\vecu_1+\vecu_2)+1
 -\epsilon_{\rho_+(\veca)}+\epsilon_{\veca}
 -\epsilon_{\rho_+(\veca+\vecu_1+\vecu_2)}
 +\epsilon_{\veca+\vecu_1+\vecu_2}
 +\epsilon_{\rho_+(\vecu_1)}\cdot\epsilon_{\rho_+(\vecu_2)}\cdot
   \epsilon_{\rho_+(\veca)}\\
=
 d(\vecv)+\rho_0(\vecu_1+\vecu_2)+2
 +
 \Bigl(
 -|\veca|_+-1
 -\epsilon_{\rho_+(\veca)}+\epsilon_{\veca}
 -\epsilon_{\rho_+(\veca+\vecu_1+\vecu_2)}
 +\epsilon_{\veca+\vecu_1+\vecu_2}
 +\epsilon_{\rho_+(\vecu_1)}\cdot\epsilon_{\rho_+(\vecu_2)}\cdot
   \epsilon_{\rho_+(\veca)}
 \Bigr).
\end{multline}
Hence we are done.
\hfill\qed

\vspace{.1in}
We need the following estimate.
\begin{lem}
Let $\vecv,\vecu_i$ and $\veca$ be as above.
We have the following inequality:
\begin{equation} \label{eq;7.4.2}
 -|\veca_i|_+ -1-\epsilon_{\rho_+(\veca)}+\epsilon_{\veca}
 -\epsilon_{\rho_+(\veca+\vecu_1+\vecu_2)}+\epsilon_{\veca+\vecu_1+\vecu_2}
 +\epsilon_{\rho_+(\vecu_1)}\cdot\epsilon_{\rho_+(\vecu_2)}
   \cdot\epsilon_{\rho_+(\veca)}
\leq g(\vecu_1)+g(\vecu_2).
\end{equation}
It is an equality if and only if
$\vecu_1\in C_L$ and $\vecu_1+\vecu_2=\veca=0$.
\end{lem}
\pf
We divide the claim into the five cases:
(1) $\vecu_i\in C_L $ $(i=1,2)$,\,\,
(2) $\vecu_i\in -C_L$ $(i=1,2)$,\,\,
(3) $\vecu_1\in -C_L$ and $\vecu_2\in C_L$,\,\,
(4) $\vecu_1\in C_L$, $\vecu_2\in (-C_L)$ and $\rho_+(\vecu_2)=0$,\,\,
(5) $\vecu_1\in C_L$, $\vecu_2\in (-C_L)$ and $\rho_+(\vecu_2)\neq 0$.

\vspace{.1in}
\noindent
{\bf Case (1)}
We have the equality $g(\vecu_1)+g(\vecu_2)=0$.
Since $\veca+\vecu_1+\vecu_2\neq 0$,
the left hand side of (\ref{eq;7.4.2})
is equal to the following:
\[
 -|\veca|_+ -1-\epsilon_{\rho_+(\veca)}+\epsilon_{\veca}
 -\epsilon_{\rho_+(\veca+\vecu_1+\vecu_2)}
 +\epsilon_{\rho_+(\vecu_1)}\cdot\epsilon_{\rho_+(\vecu_2)}
   \cdot\epsilon_{\rho_+(\veca)}.
\]
We have the inequalities
$-\epsilon_{\rho_+(\veca)}+\epsilon_{\veca}\leq 0 $
and $-\epsilon_{\rho_+(\veca+\vecu_1+\vecu_2)}
 +\epsilon_{\rho_+(\vecu_1)}\cdot\epsilon_{\rho_+(\vecu_2)}
   \cdot\epsilon_{\rho_+(\veca)}\leq 0 $.
Thus the left hand side of (\ref{eq;7.4.2}) is less than $-1$.

\noindent
{\bf Case (2)}
Since $\veca+\vecu_1+\vecu_2$ is contained in $\nbiga_L$,
we know $\veca\neq 0$ and 
$-|\veca|_+\leq g(\vecu_1)+g(\vecu_2)$.
Thus the left hand side of (\ref{eq;7.4.2}) is less than the following:
\[
 g(\vecu_1)+g(\vecu_2)-1
 -\epsilon_{\rho_+(\veca)}-\epsilon_{\rho_+(\veca+\vecu_1+\vecu_2)}
 +\epsilon_{\veca+\vecu_1+\vecu_2}
 +\epsilon_{\rho_+(\vecu_1)}\cdot\epsilon_{\rho_+(\vecu_2)}
  \cdot\epsilon_{\rho_+(\veca)}.
\]
Since we have the inequalities
$-\epsilon_{\rho_+(\veca+\vecu_1+\vecu_2)}
 +\epsilon_{\veca+\vecu_1+\vecu_2}\leq 0$
and
$-\epsilon_{\rho_+(\veca)}+\epsilon_{\rho_+(\vecu_1)}\cdot\epsilon_{\rho_+(\vecu_2)}
  \cdot\epsilon_{\rho_+(\veca)}\leq 0$,
the left hand side of (\ref{eq;7.4.2})
is strictly less than $g(\vecu_1)+g(\vecu_2)$.

\noindent
{\bf Case (3)} 
Since $\vecu_1+\veca$ is contained in $\nbiga_L$,
we know that $\veca\neq 0$
and that $-|\veca|_+\leq g(\vecu_1)=g(\vecu_1)+g(\vecu_2)$.
Thus the left hand side of (\ref{eq;7.4.2})
is less than the following:
\[
 g(\vecu_1)+g(\vecu_2)-1
 -\epsilon_{\rho_+(\veca)}-\epsilon_{\rho_+(\veca+\vecu_1+\vecu_2)}
 +\epsilon_{\veca+\vecu_1+\vecu_2}
 +\epsilon_{\rho_+(\vecu_1)}\cdot\epsilon_{\rho_+(\vecu_2)}
  \cdot\epsilon_{\rho_+(\veca)}.
\]
As in the case (2), it is strictly less than $g(\vecu_1)+g(\vecu_2)$.

\noindent
{\bf Case (4)}
We have $g(\vecu_1)+g(\vecu_2)=0$.
If $\veca\neq 0$,
the left hand side of (\ref{eq;7.4.2}) is less than the following:
\begin{equation} \label{eq;7.8.1}
 -|\veca|_+ -1-\epsilon_{\rho_+(\veca)}
 -\epsilon_{\rho_+(\veca+\vecu_1+\vecu_2)}
 +\epsilon_{\veca+\vecu_1+\vecu_2}
 +\epsilon_{\rho_+(\vecu_1)}\cdot\epsilon_{\rho_+(\veca)}.
\end{equation}
Since we have
$-\epsilon_{\rho_+(\veca)}+
 \epsilon_{\rho_+(\vecu_1)}\cdot \epsilon_{\rho_+(\veca)}\leq 0$
and
$-\epsilon_{\rho_+(\veca+\vecu_1+\vecu_2)}
 +\epsilon_{\veca+\vecu_1+\vecu_2}\leq 0$,
the left hand side of (\ref{eq;7.4.2}) is less than $-1$.

If $\veca=0$,
the left hand side of (\ref{eq;7.4.2}) is equal to the following:
\[
 0-1-1+1-\epsilon_{\rho_+(\vecu_1+\vecu_2)}
        +\epsilon_{\vecu_1+\vecu_2}
        +\epsilon_{\rho_+(\vecu_1)}=-1+\epsilon_{\vecu_1+\vecu_2}.
\]
It is less than $0=g(\vecu_1)+g(\vecu_2)$,
and it is an equality if and only if $\vecu_1+\vecu_2=0$.

\noindent
{\bf Case (5)}
We have the equality $g(\vecu_1)+g(\vecu_2)=-1$.
If $\veca\neq 0$,
then the left hand side of (\ref{eq;7.4.2}) is less than
the following:
\[
 -|\veca|_+ -1-\epsilon_{\rho_+(\veca)}
 -\epsilon_{\rho_+(\veca+\vecu_1+\vecu_2)}+\epsilon_{\veca+\vecu_1+\vecu_2}
 \leq
 -|\veca|_+ -1-\epsilon_{\rho_+(\veca)}
 \leq -2<g(\vecu_1)+g(\vecu_2).
\]
Assume that $\veca=0$,
the left hand side is equal to the following:
\[
 0-1-1+1-\epsilon_{\rho_+(\vecu_1+\vecu_2)}+\epsilon_{\vecu_1+\vecu_2}+0
 =-1-\epsilon_{\rho_+(\vecu_1+\vecu_2)}+\epsilon_{\vecu_1+\vecu_2}.
\]
Our assumption
that $\veca+\vecu_1+\vecu_2=\vecu_1+\vecu_2\in \nbiga_L$
implies the inequality
$\rho_{\alpha}(\vecu_{1}+\vecu_{2})\geq 0$ for any $\alpha>0$.
Since we have $\vecu_1\in C_L$ and $\vecu_2\in -C_L$,
we obtain the equality $\rho_{\alpha}(\vecu_{1}+\vecu_{2})=0$
for any $\alpha>0$,
i.e.,
$\rho_+(\vecu_1+\vecu_2)=0$.
Thus the left hand side of (\ref{eq;7.4.2}) is less than
$-2+\epsilon_{\vecu_1+\vecu_2}$.
It is less than $-1=g(\vecu_1)+g(\vecu_2)$.
It is equal to $-1$ if and only if $\vecu_1+\vecu_2=0$.
\hfill\qed

We obtain the following direct corollary.
\begin{cor} \label{cor;3.10.2}
Let $(\vecv,\vecu_1,\vecu_2)$ be an elements
of $Comp_L$.
Let $\veca$ be an element of $St_1(\vecv,\vecu_1,\vecu_2)$.
Then we have the following inequality:
\[
 \dim(\widetilde{Z}(\vecv,\vecu_1,\vecu_2)_{1,\veca})
\leq
 d(\vecv,\vecu_1,\vecu_2).
\]
It is an equality if and only if 
$\vecu_1\in C_L$ and $\vecu_1+\vecu_2=\veca=0$.
\hfill\qed
\end{cor}

\noindent
{\bf Notation:}
When $\vecu_1\in C_L$ and $\vecu_1+\vecu_2=0$,
we denote
$\widetilde{Z}(\vecv,\vecu_1,\vecu_2)_{1,0}$
by $\Delta'$.
The closure $\overline{Z}(\vecv,\vecu_1,\vecu_2)_{1,0}$
is same as  the diagonal of
$X^{[\vecv]}\times Y_{\vecu_1}\times Y_{\vecu_1}\times X^{[\vecv]}$.
We denote it by $\Delta_{X^{[\vecv]}\times Y_{\vecu_1}}$.


\subsection{The multiplicities}
Due to Proposition \ref{prop;3.10.1},
the cycle
$\pi_{3,\ast}\bigl([Z(\vecv,\vecu_1)]\cdot[Z(\vecv+\vecu_1,\vecu_2)]\bigr)$
contained in
$CH_{d(\vecv,\vecu_1,\vecu_2)}\bigl(\overline{Z}(\vecv,\vecu_1,\vecu_2)\bigr)$
is described as follows:
\[
 \pi_{3,\ast}\bigl([Z(\vecv,\vecu_1)]\cdot[Z(\vecv+\vecu_1,\vecu_2)]\bigr)
=
 \alpha\cdot [\overline{Z}(\vecv,\vecu_1,\vecu_2)_0]+
 \beta\cdot
 h(\vecu_1,\vecu_2)\cdot \Delta_{X^{[\vecv]}\times Y_{\vecu_1}}.
\]
Here $\alpha$ and $\beta$ denote the multiplicities,
and $h(\vecu_1,\vecu_2)$ is defined to be $1$
$(\vecu_1\in C,\vecu_1+\vecu_2=0)$
or $0$ (otherwise).
We need to calculate the multiplicities.
We put as follows:
\[
W=
\left\{
\begin{array}{ll}
 \overline{Z}(\vecv,\vecu_1,\vecu_2)-
 \bigl(\widetilde{Z}(\vecv,\vecu_1,\vecu_2)_0\cup
  \Delta'
 \bigr),
 & (\vecu_1\in C,\vecu_1+\vecu_2=0),\\
 \overline{Z}(\vecv,\vecu_1,\vecu_2)-\widetilde{Z}(\vecv,\vecu_1,\vecu_2)_0,
 & (\mbox{otherwise}).
\end{array}
\right.
\]
Note that $W$ is closed.
We already know that
$\dim(W)<d(\vecv,\vecu_1,\vecu_2)$.
We put
$U:=X^{[\vecv]}\times Y_{\vecu_1}\times Y_{\vecu_2}
 \times X^{[\vecv+\vecu_1+\vecu_2]}-W$.
It is easy to see that
$U\cap \overline{Z}(\vecv,\vecu_1,\vecu_2)=
\widetilde{Z}(\vecv,\vecu_1,\vecu_2)_0\sqcup \Delta'$.
The intersection $Z(\vecv,\vecu_1,\vecu_2)\cap \pi_3^{-1}(U)$
is divided into the disjoint union
of $\pi_3^{-1}(\widetilde{Z}(\vecv,\vecu_1,\vecu_2)_0)$
and $\pi_3^{-1}(\Delta')$.

The subvarieties $\pi_{4,5}(Z(\vecv,\vecu_1))$
and $\pi_{1,2}(Z(\vecv+\vecu_1,\vecu_2))$
intersect transversally at
$\pi_3^{-1}(\widetilde{Z}(\vecv,\vecu_1,\vecu_2)_0)$.
Thus the multiplicity $\alpha$ is $1$.

Let consider the multiplicity $\beta$.
Let $\omega$ be the natural morphism
$X^{[\vecu_1]}\lrarr X^{(\vecu_1)}$.
We have the natural diagonal embedding
$Y_{\vecu_1}\lrarr X^{(\vecu_1)}$.
Let take a point $x$ of $Y_{\vecu_1}$,
and then
it gives the point of $X^{(\vecu_1)}$.
We put $X^{[\vecu_1]}_0:=\omega^{-1}(Y_{\vecu_1})$
and $X^{[\vecu_1]}_x:=\omega^{-1}(x)\subset X^{[\vecu_1]}_0$.

The number $\beta$ is same as the intersection number
$[X^{[\vecu_1]}_{x}]\cdot [X^{[\vecu_1]}_0]$
in $X^{[\vecu_1]}$.
When $\rho_+(\vecu_1)=0$,
the number $[X^{[\vecu_1]}_{x}]\cdot [X^{[\vecu_1]}_0]$
is calculated by Ellingsrud and Str{\o}mme,
which is $(-1)^{\rho_0(\vecu_1)}\cdot \rho_0(\vecu_1)$.
The following lemma can be easily derived from their calculation.
\begin{lem}
When $\rho_+(\vecu_1)\neq 0$,
the number $[X^{[\vecu_1]}_{x}]\cdot [X^{[\vecu_1]}_0]$
is $(-1)^{\rho_0(\vecu_1)}$.
\end{lem}
\pf
We put $\rho_0(\vecu_1)=n$.
Let $X^{[n,n+1]}$ denote the nested Hilbert scheme,
that is,
the moduli of the pairs
$(I_1,I_2)$ of the ideal sheaves on $X$
such that $\length(\nbigo/I_1)=n,\length(\nbigo/I_2)=n+1$
and that $I_1\supset I_2$.
It is well known that $X^{[n,n+1]}$ is smooth.
We have the morphism $\omega:X^{[n,n+1]}\lrarr X^{(n+1)}$.
We denote the small diagonal of $X^{(n+1)}$
by $\Delta$.
We denote $\omega^{-1}(\Delta)$ by $X^{[n,n+1]}_0$.
Let take a point $x$ of $\Delta$.
We denote $\omega^{-1}(x)$ by $X^{[n,n+1]}_x$.
Ellingsrud-Str{\o}mme calculated the intersection number
$[X^{[n,n+1]}_0]\cdot[X^{[n,n+1]}_x]$
in $X^{[n,n+1]}$.
The answer is $(-1)^n$.

We have the natural morphism $H:X^{[n,n+1]}\lrarr X$.
The morphism $H$ is transversal with $D$.
Moreover the restriction of $H$
to $X^{[n,n+1]}_0$ is also transversal with $D$.

If $\vecu$ is an element of $C$
such that $\rho_+(\vecu)\neq 0$,
then $X^{[\vecu]}$ (resp. $X^{[\vecu]}_0$)
is naturally isomorphic to $H^{-1}(D)$
(resp. $H|_{X^{[n,n+1]}_0}^{-1}(D)$).
Thus the intersection number desired is $(-1)^n$.
\hfill\qed

In all, we obtain the following result.
\begin{thm} \label{thm;3.3.5}
Let 
$(\vecv,\vecu_1,\vecu_2)$ be an element of $Comp_L$.
In the group $CH_{d(\vecv,\vecu_1,\vecu_2)}(\pi_3(Z(\vecv,\vecu_1,\vecu_2)))$,
we have the following equality:
\[
 \pi_{3\,\ast}\bigl(
 [Z(\vecv,\vecu_1)]\cdot[Z(\vecv+\vecu_1,\vecu_2)]
 \bigr)=
 [\overline{Z}(\vecv,\vecu_1,\vecu_2)_0]+
  f(\vecu_1,\vecu_2)\cdot
 \bigl(-\rho_0(\vecu_1)\bigr)^{\epsilon_{\rho_+(\vecu_1)}}\cdot
 (-1)^{\rho_0(\vecu_1)}\cdot
 \Delta_{X^{[\vecv]}\times Y_{\vecu_1}}.
\]
Here $f(\vecu_1,\vecu_2)$ is defined to be $1$
$(\vecu_1\in C,\vecu_1+\vecu_2=0)$
or $0$ (otherwise).
\hfill\qed
\end{thm}

\section{The Heisenberg relation}

We consider the parabolic Hilbert schemes $X^{[\vecv]}$
such that $\vecv$ is not necessarily contained in $\nbiga_L$.

\subsection{The composition of incidence varieties}

Let $\vecv$ be an element of $\nbiga$.
We put
$\tilde{C}(\vecv):=\{\vecu\in\tilde{C}\,|\,\vecv+\vecu\in \nbiga\}$
and $\vecu$ be an element of $\tilde{C}(\vecv)$.
Then we have already introduced the incidence variety
$Z(\vecv,\vecu)$ for $(\vecv,\vecu)$
in subsubsection \ref{subsub;7.8.10}.
The following lemma can be shown by an argument similar
to the proof of Lemma \ref{lem;3.1.1}
and Lemma \ref{lem;7.8.15}.
\begin{lem}
The dimension of the incidence variety $Z(\vecv,\vecu)$
is $d(\vecv)+\rho_0(\vecu)+1+g(\vecu)$.
Here the function $g:\tilde{C}\lrarr \{-1,1\}$
is defined in subsubsection {\rm \ref{subsubsection;7.6.10}}.
\hfill\qed
\end{lem}

For $\vecv$,
the subset $\widetilde{C}^2(\vecv)$ of $\widetilde{C}^2$
is defined to be
$\{(\vecu_1,\vecu_2)\in \widetilde{C}^2\,|\,
 \vecv+\vecu_1\in \nbiga,\,\vecv+\vecu_1+\vecu_2\in \nbiga\}$.
The subset $Comp$ of $\nbiga\times \widetilde{C}^2$
is defined to be
$\{(\vecv,\vecu_1,\vecu_2)\,|\,\vecv\in \nbiga,\,(\vecu_1,\vecu_2)\in
 \widetilde{C}^2(\vecv)\}$.
For any element $(\vecv,\vecu_1,\vecu_2)\in Comp$,
we have the closed subscheme
$\pi_3(Z(\vecv,\vecu_1,\vecu_2))$
of $X^{[\vecv]}\times Y_{\vecu_1}\times Y_{\vecu_2}\times
 X^{[\vecv+\vecu_1+\vecu_2]}$
as in subsection \ref{subsection;7.8.20}.
Moreover
the closed subscheme
$\overline{Z}(\vecv,\vecu_1,\vecu_2)_0$ of
$X^{[\vecv]}\times Y_{\vecu_1}
 \times Y_{\vecu_2}\times X^{[\vecv+\vecu_1+\vecu_2]}$
is defined to be the closure of the following locally closed subscheme:
\[
 \left\{
 \bigl((I_1,\nbigg_{1,\ast}),x_1,x_2,
       (I_2,\nbigg_{2,\ast})
 \bigr)\,\Bigl|\Bigr.\,
\begin{array}{l}
 x_1\neq x_2,\,\,
 (I_1,\nbigg_{1,\ast})=(I_2,\nbigg_{2,\ast})
 \mbox{ on $X-\{x_1,x_2\}$, }\\
 \length_{x_i}(\nbigo/I_1,\nbigg_{1\,\ast})=\max(0,-\vecu_i),\,
 \length_{x_i}(\nbigo/I_2,\nbigg_{2\,\ast})=\max(0,\vecu_i)
 \end{array}
 \right\}.
\]
Here $\max(0,\vecu)$ denotes $\vecu$ $(\vecu\in C)$
or $0$ $(\vecu\in -C)$ for any $\vecu\in \tilde{C}$.
If $(\vecv,\vecu_1,\vecu_2)$ is contained in
$Comp_L$,
then $\overline{Z}(\vecv,\vecu_1,\vecu_2)_0$
is same as the variety considered in the previous section.
It is easy to see that the composition varieties
are preserved by the morphisms of the form $\Psi_{\vecv,\vecv'(\beta)}$
given in subsection \ref{subsection;7.4.10}.
If  $\beta=-1$,
then the element $\vecv'(-1)$ associated to
$\vecv\in \nbiga\cap N(\alpha_-,\alpha_+)$
is an element of $\nbiga_L$.
The following theorem is a direct corollary of Theorem \ref{thm;3.3.5}.
\begin{thm}
Let 
$(\vecv,\vecu_2,\vecu_2)$ be an element of $Comp$.
We put
$d(\vecv,\vecu_1,\vecu_2):=
 d(\vecv)+\rho_0(\vecu_1+\vecu_2)+2+g(\vecu_1)+g(\vecu_2)$.
Then $\pi_3(Z(\vecv,\vecu_1,\vecu_2))$ is
$d(\vecv,\vecu_1,\vecu_2)$-dimensional.
In the Chow group
$CH_{d(\vecv,\vecu_1,\vecu_2)}(
 \pi_3(Z(\vecv,\vecu_1,\vecu_2)))$,
we have the following equality:
\begin{multline}
 \pi_{3\,\ast}\bigl(
 [Z(\vecv,\vecu_1)]\cdot[Z(\vecv+\vecu_1,\vecu_2)]
 \bigr)
\\
=[\overline{Z}(\vecv,\vecu_1,\vecu_2)_0]+
f (\vecu_1,\vecu_2)\cdot
 \bigl(-\rho_0(\vecu_1)+|\vecu_1|_-\bigr)^{\epsilon_{\rho(\vecu_1)}}\cdot
 (-1)^{\rho_0(\vecu_1)-|\vecu_1|_-}\cdot
 \Delta_{X^{[\vecv]}\times Y_{\vecu_1}}.
\end{multline}
Here $f(\vecu_1,\vecu_2)$ is defined to be $1$
$(\vecu_1\in C,\,\vecu_1+\vecu_2=0)$
or $0$ (otherwise).
\hfill\qed
\end{thm}

\subsection{The commutator of the compositions}

Formally we can consider the parabolic Hilbert schemes
for any $\vecv\in \TT$:
Since there exists no parabolic ideal
with parabolic length $\vecv\in \TT-\nbiga$,
the variety $X^{[\vecv]}$ is defined to be
the empty set for such $\vecv$.
Let $\vecv$ be an element of $\TT$
and $\vecu$ be an element of $\widetilde{C}$.
When one of $\vecv$ or $\vecv+\vecu$
is not contained in $\nbiga$,
the incidence variety
$Z(\vecv,\vecu)$ is defined to be the empty set.
Then we naturally obtain the varieties
$\pi_3(Z(\vecv,\vecu_1,\vecu_2))$
and $\overline{Z}(\vecv,\vecu_1,\vecu_2)_0$
for any $\vecv\in \TT$ and $(\vecu_1,\vecu_2)\in \widetilde{C}^2$.

Let $(\vecv,\vecu_1,\vecu_2)$
be any element of $\TT\times \widetilde{C}^2$.
We have the two varieties
$X^{[\vecv]}\times Y_{\vecu_1}\times Y_{\vecu_2}
 \times X^{[\vecv+\vecu_1+\vecu_2]}$
and
$X^{[\vecv]}\times Y_{\vecu_2}\times Y_{\vecu_1}
 \times X^{[\vecv+\vecu_1+\vecu_2]}$.
We denote the transposition of the second and third component
by $\tau_{2\,3}$:
\[
 \tau_{2\,3}:
 X^{[\vecv]}\times Y_{\vecu_1}\times Y_{\vecu_2}
 \times X^{[\vecv+\vecu_1+\vecu_2]}
\lrarr
 X^{[\vecv]}\times Y_{\vecu_2}\times Y_{\vecu_1}
 \times X^{[\vecv+\vecu_1+\vecu_2]},
\quad
(x_1,x_2,x_3,x_4)\longmapsto (x_1,x_3,x_2,x_4).
\]
Let $\pi_3$ denote the projection of
$X^{[\vecv]}\times Y_{\vecu_1}\times X^{[\vecv+\vecu_1]}
\times Y_{\vecu_2}\times X^{[\vecv+\vecu_1+\vecu_2]}$
and
$X^{[\vecv]}\times Y_{\vecu_2}\times X^{[\vecv+\vecu_2]}
\times Y_{\vecu_1}\times X^{[\vecv+\vecu_1+\vecu_2]}$
onto
$X^{[\vecv]}\times Y_{\vecu_1}
\times Y_{\vecu_2}\times X^{[\vecv+\vecu_1+\vecu_2]}$
and 
$X^{[\vecv]}\times Y_{\vecu_2}
\times Y_{\vecu_1}\times X^{[\vecv+\vecu_1+\vecu_2]}$
respectively.
For each element $\vecu\in \tilde{C}$,
the number $sgn(\vecu)$ is defined to be $1$ ($\vecu\in C$)
or $-1$ ($\vecu\in -C$).

\begin{thm} \label{thm;3.3.1}
Let $\vecv$ be an element of $\nbiga$ and $(\vecu_1,\vecu_2)$
be an element of $\widetilde{C}^2$.
Then we have the following relation
in the group
$CH_{d(\vecv,\vecu_1,\vecu_2)}
   (X^{[\vecv]}\times Y_{\vecu_1}\times Y_{\vecu_2}
    \times X^{[\vecv+\vecu_1+\vecu_2]})$:
\begin{multline}
 \pi_{3\,\ast}\bigl([Z(\vecv,\vecu_1)]\cdot [Z(\vecv+\vecu_1,\vecu_2)]\bigr)
-\tau_{2\,3\,\ast}\pi_{3\,\ast}
 \bigl([Z(\vecv,\vecu_2)]\cdot [Z(\vecv+\vecu_2,\vecu_1)]\bigr)\\
=\epsilon_{\vecu_1+\vecu_2}\cdot
 sgn(\vecu_1)
 \cdot
 (-|\rho_0(\vecu_1)|+|\vecu_1|_- )^{\epsilon_{\rho(\vecu_1)}}
\cdot
 (-1)^{\rho_0(\vecu_1)-|\vecu_1|_-}\cdot
 \Delta_{X^{[\vecv]}\times Y_{\vecu_1}}.
\end{multline}
\end{thm}
\pf
It is easy to see that
$\widetilde{Z}(\vecv,\vecu_1,\vecv_2)_0$
and $\tau_{2\,3}\bigl(\widetilde{Z}(\vecv,\vecu_2,\vecu_1)_0\bigr)$
are same
for any $(\vecv,\vecu_1,\vecu_2)$
by definition.
Thus their closures are same.
Hence the terms
$[\overline{Z}(\vecv,\vecu_1,\vecu_2)_0]$
and the
$\tau_{2\,3\,\ast}\bigl([\widetilde{Z}(\vecv,\vecu_2,\vecu_1)_0]\bigr)$
are canceled.
\hfill\qed


Let $\vecu$ be an element of $\tilde{C}$ and
$a$ be an element of $H^{\ast}(Y_{\vecu})$.
Let consider the operator $\gminiq_{\vecu}(a)$
defined in subsection \ref{subsubsection;3.3.3}.
We obtain the following direct corollary.
\begin{cor}  \label{cor;3.10.5}
Let $\vecu_1,\vecu_2$ be elements of $\tilde{C}$.
Let $a_i$ be an element of $H^{\ast}(Y_{\vecu_i})$
for $i=1,2$.
Then the following relation holds
\[
 \bigl[\gminiq_{\vecu_1}(a_1),\gminiq_{\vecu_2}(a_2)\bigr]
=
 \epsilon_{\vecu_1+\vecu_2}\cdot
 \mu(\vecu_1)
 \cdot
 \left(
 \int_{Y_{\vecu_1}}a_1\cdot a_2\right)
 \cdot
 id_{\hyperh(X,D)}.
\]
Here $id_{\hyperh(X,D)}$ denotes
the identity of $\hyperh(X,D)$,
and we put as follows:
\[
 \mu(\vecu_1)=
 sgn(\vecu_1)\cdot
 \bigl(-|\rho_0(\vecu_1)|+|\vecu_1|_- \bigr)
  ^{\epsilon_{\rho(\vecu_1)}}
 (-1)^{\rho_0(\vecu_1)-|\vecu_1|_-}.
\]
In particular, the number $\mu(\vecu_1)$ is not $0$.
\hfill\qed
\end{cor}



\noindent
{\it Address\\
Department of Mathematics,
Osaka City University,
Sugimoto, Sumiyoshi-ku,
Osaka 558-8585, Japan.
takuro@sci.osaka-cu.ac.jp\\
School of Mathematics,
Institute for Advanced Study,
Einstein Drive, Princeton,
08540, NJ, USA. \\
takuro@math.ias.edu
}

\end{document}